\documentclass[12pt,a4paper,draft]{article}

\usepackage{amsmath,amsfonts}
\usepackage{mathrsfs} 
\usepackage[T1]{fontenc}
\usepackage[utf8]{inputenc}
\usepackage{authblk}
\usepackage{enumerate}
\usepackage{graphicx}
\usepackage{xcolor}

\newtheorem{theorem}{Theorem}
\newtheorem{lemma}{Lemma}

\let\scr\mathscr

\def \Liminf{\mathop{\underline{\lim}}\limits}

\def\1{\mbox{1\hspace{-.25em}I}}

\def\Pb{\mathbf{P}}

\def\Ex{\mathbf{E}}

\def\UU{\mathbb{U}}

\def\KK{\mathbb{K}}

\def\1{\mbox{1\hspace{-.25em}I}}

\begin{document}

\title{On Multi-step Estimation of  Delay for SDE  } 
\author{\textsc{Yury A. Kutoyants}\\ 
{\small Le Mans University, Le Mans,    France }\\
{\small Tomsk     State University, Tomsk, Russia }
}

\date{}

\maketitle

\begin{abstract}
We consider the problem of delay estimation by the observations of the
solutions of several SDEs. It is known that the MLE for these models are
consistent and asymptotically normal, but the likelihood ratio functions are
not differentiable w.r.t. parameter and therefore numerical calculation of the
MLEs has certain difficulties. We propose One-step and Two-step MLE, whose
calculation has no such problems and provide estimator asymptotically
equivalent to the MLE. These constructions are realized in two or three
steps. First we construct preliminary estimators which are consistent and
asymptotically normal, but not asymptotically efficient.   Then we use these
estimators and modified Fisher-score device to obtain One-step and Two-step
MLEs. We suppose that its numerical realization is much more
simple. Stochastic Pantograph equation is introduced and related statistical
problems are discussed.
\end{abstract}

\section{Introduction}

We consider two types of models of stochastic differential equation (SDE) with
delays, which can be illustrated by the following equations
\begin{align}
\label{1}
{\rm d}X_t&=S\left(X_{t-\vartheta }\right){\rm d}t +\varepsilon {\rm d}W_t, \quad
X_s=x_0, \;s\leq 0, \qquad 0\leq t\leq T,\\
{\rm d}X_t&=S\left(X_{\vartheta t}\right){\rm d}t +\varepsilon {\rm d}W_t, \qquad
X_0=x_0, \;, \qquad 0\leq t\leq T.
\label{1a}
\end{align}
We call delay in the equation \eqref{1} {\it shift-type} and in the equation
\eqref{1a} {\it scale-type}. 
 Here the functions $S\left(x\right)$ are
known and smooth and parameter (delay) $\vartheta \in \Theta $.  The set
$\Theta=\left(\alpha ,\beta \right)$, $0<\alpha< \beta <T$ in the case of
equation \eqref{1} and $\Theta=\left(\alpha ,\beta \right)$, $0<\alpha< \beta
<1$ in the case of equation \eqref{1a}. The parameter $\vartheta $ has to be
estimated 
by observations $X^T=\left(X_t,0\leq t\leq T\right)$ and the properties of
estimators are studied in the asymptotic of {\it small noise}, i.e., as $\varepsilon \rightarrow 0 $. Such perturbed dynamical
systems were studied in \cite{FW84} and different statistical problems related
to such models can be found in \cite{Kut94}. We suppose that the function
$S\left(x\right)$ is Lipschitz 
\begin{align*}
\left|S\left(x\right)-S\left(y\right)\right|\leq L\left|x-y\right|.
\end{align*}
Under this condition the equations \eqref{1} and \eqref{1a} have unique strong
solutions and all polynomial moments of the solutions are finite (see Theorem
4.6 in \cite{LS01}).

Ordinary and stochastic differential equations with delay (feedback with
delay) are widely used in many applied sciences (see \cite{BLMSS05} (medicine),
\cite{FFB05}, \cite{OY00},  \cite{RTM18} (physics), \cite{HPT05} (neuron science),
\cite{Ku93}, \cite{M05},  (population
dynamics)
and references therein). The same time statistical inference (estimation and
testing concerning delays) for these systems does not attract sufficiently
attention of the statisticians.

The possible difficulties in the study of these models of observations we
illustrate with the help of the linear Ornstein-Uhlenbeck process with delay
\begin{align}
\label{2}
{\rm d}X_t=-\gamma X_{t-\vartheta  }{\rm d}t +\sigma  {\rm d}W_t, \quad
X_s=x_0, \;s\leq 0,\quad 0\leq t\leq T.  
\end{align}
Consider the problem of estimation $\vartheta \in\Theta$ by continuous time
observations $X^T=\left(X_t,0\leq t\leq T\right)$ from \eqref{2}. We suppose
that the values $\gamma >0$ and $\sigma >0$ are known.

 Let us write the equation \eqref{2} in integral form and
change the variables $t\rightarrow t-\vartheta $. Then we obtain the
representation
\begin{align}
\label{3}
X_{t-\vartheta }=x_0 -\gamma \int_{0}^{t-\vartheta }X_{s-\vartheta }{\rm
  d}s+\sigma W_{t-\vartheta },\qquad t\geq \vartheta . 
\end{align}
We see that the drift $S\left(X_{t-\vartheta })\right)=-\gamma X_{t-\vartheta }$
in the equation \eqref{2} is as smooth w.r.t. $\vartheta $ as the Wiener
process w.r.t. time $t$, i.e., even the first derivative of the drift
$\partial S\left(X_{t-\vartheta })\right)/\partial \vartheta $ does not
exist. The following question naturally arises: does the problem of
estimation $\vartheta $ by observations \eqref{3} {\it regular} (family of
measures is LAN, Fisher information is finite), or {\it singular} (family of
measures is not LAN, Fisher information is infinite)? What is amusing with
this model of observations, we have the both cases depending on the type of
asymptotics. If the time of observations is fixed and $\sigma \rightarrow 0$,
then the problem of estimation is {\it regular}, family of measures is LAN and
MLE $\hat\vartheta _\sigma $ is asymptotically normal
\begin{align*}
\sigma ^{-1}\left(\hat\vartheta _\sigma-\vartheta \right)\Longrightarrow {\cal
N}\left(0,{\rm I}\left(\vartheta \right)^{-1}\right),\qquad {\rm
  I}\left(\vartheta \right) =\gamma ^4\int_{\vartheta }^{T}x_{t-2\vartheta }^2{\rm d}t
\end{align*}
(see \cite{Kut88}, \cite{Kut94}). Here ${\rm I}\left(\vartheta \right)$ plays
the role of Fisher information and $x_t=x_t\left(\vartheta \right)$ is
solution of the limit ($\sigma =0$) equation
\begin{align}
\label{4}
\frac{{\rm d}x_t}{{\rm d}t}=-\gamma x_{t-\vartheta },\qquad x_s=x_0,\;s\leq
0,\qquad 0\leq t\leq T.
\end{align}

 Note that the observations $X_t$ on the
interval $\left[0,\vartheta \right]$ has no any information about the value of
$\vartheta $ that is why the derivative of $-\gamma X_{t-\vartheta }=-\gamma
x_0$ w.r.t. $\vartheta $ is zero.  

  If we fix
$\sigma $ and consider asymptotic $T\rightarrow \infty $, then the problem of
estimation becomes {\it singular}, family of measures is no more LAN and MLE
$\hat\vartheta _T$ has another (non Gaussian) limit distribution
\begin{align*}
Tc\left(\hat\vartheta _T-\vartheta \right)\Longrightarrow
{\rm arg}\sup_{u\in {\cal R}}\left(W\left(u\right)-\frac{\left|u\right|}{2}\right) .
\end{align*}
Here $c>0$ is some constant and $W\left(\cdot \right)$ is two-sided Wiener
process \cite{KuK00}, \cite{Kut04}. Recall that the similar properties has MLE
$\hat\vartheta _T$ in the case of observations
\begin{align*}
{\rm d}X_t=-\gamma\; {\rm sgn} \left(X_{t }-\vartheta\right){\rm d}t +\sigma
{\rm d}W_t, \quad X_s=x_0, \;s\leq 0,\qquad 0\leq t\leq T
\end{align*}
with discontinuous drift coefficient \cite{Kut04}. Hence here delay is
equivalent to discontinuity (change-point in space) in the model. This is the
only known to us model, which is the same time regular and singular
(statistically) depending on the type of limit. The proofs are based on the
powerful general results in \cite{IH81}.

The properties of the MLE $\hat\vartheta _\varepsilon $ for the model
\eqref{1} were described in the works \cite{Kut88} (linear case) and in
\cite{A86} (generalization of \cite{Kut88} to nonlinear case), see as well
\cite{Kut94}. It was shown that the MLE is consistent, asymptotically normal
\begin{align*}
\varepsilon ^{-1}\left(\hat\vartheta _\varepsilon -\vartheta \right)\Longrightarrow {\cal
N}\left(0,{\rm I}\left(\vartheta \right)^{-1}\right),\qquad {\rm
  I}\left(\vartheta \right) =\int_{\vartheta }^{T}S'\left(x_{t-\vartheta
}\right)^2S\left(x_{t-2\vartheta }\right)^2{\rm d}t 
\end{align*} 
and asymptotically efficient.

The likelihood ratio function for the model \eqref{1} is   
\begin{align}  
\label{LR}
L\left(\vartheta ,X^T\right)=\exp\left\{\int_{0}^{T}\frac{S\left(X_{t-\vartheta
}\right)}{\varepsilon ^2}{\rm d}X_t-\int_{0}^{T}\frac{S\left(X_{t-\vartheta
}\right)^2}{2\varepsilon ^2}{\rm d}t\right\},\quad \vartheta \in\Theta 
\end{align}
and the MLE $\hat\vartheta _\varepsilon $  is defined by the equation 
\begin{align}
\label{mle}
L(\hat\vartheta _\varepsilon  ,X^T)=\sup_{\vartheta \in\Theta }L\left(\vartheta ,X^T\right).
\end{align}
As the function $L\left(\vartheta ,X^T\right)$ is not differentiable
w.r.t. $\theta $ the realization of numerical algorithms of calculation $\hat\vartheta
_\varepsilon  $ can have certain difficulties.

The goal of this work is to construct estimators, whose calculation has no
such difficulties and which have the same asymptotic properties as the MLE for
the both models. This will be realized in two or three steps. First we
construct some preliminary estimators, which have sufficiently good rate of
convergence, but not asymptotically efficient. Then using these estimators
and Fisher-score device we obtain One-step MLEs $\vartheta ^\star_\varepsilon
$. It is shown that the estimators $\vartheta ^\star_\varepsilon $ are
consistent and have the same asymptotic variance as the MLEs. The Fisher-score
improvement  of the preliminary estimators is a well-known approach in
statistics and we do not give here the review of the corresponding
literature. The discussion of the works related with diffusion processes can
be found in \cite{Kut17}. 

Recall that the realization of Fisher-score device requires two derivatives of
the trend coefficient w.r.t. unknown parameter. For example, suppose that the
observed process is
\begin{align*}
{\rm d}X_t=S\left(\vartheta ,t\right){\rm d}t+\varepsilon  {\rm d}W_t,\quad
X_0=0,\quad 0\leq t\leq T 
\end{align*}
and we have preliminary estimator $\bar\vartheta_\varepsilon  $ such that
$\varphi _\varepsilon ^{-1}\left(\bar\vartheta_\varepsilon-\vartheta _0\right)
$ is bounded in probability. Then we can write the One-step MLE formally as
follows 
\begin{align*}
\vartheta _\varepsilon ^\star=\bar\vartheta_\varepsilon+ \int_{0}^{T} \frac{\dot
S\left(\bar\vartheta_\varepsilon,t\right)}{ {\rm
  I}\left(\bar\vartheta_\varepsilon\right)  }\left[{\rm
    d}X_t-S\left(\bar\vartheta_\varepsilon,t\right){\rm d}t\right],\qquad {\rm
  I}\left(\vartheta\right)=\int_{0}^{T} \dot
S\left(\vartheta,t\right)^2{\rm d}t.
\end{align*}
Hence, using Taylor formula
$$
S\left(\vartheta_0,t\right)-S\left(\bar\vartheta_\varepsilon,t\right)=-\left(\bar\vartheta_\varepsilon-\vartheta_0\right)\dot  
S\left(\bar\vartheta_\varepsilon,t\right)
-\frac{1}{2}\left(\bar\vartheta_\varepsilon-\vartheta_0\right)^2\ddot
S(\tilde\vartheta_\varepsilon,t)
$$
we obtain the representation
\begin{align*}
\frac{\vartheta _\varepsilon ^\star-\vartheta _0}{\varepsilon }& =\varepsilon ^{-1}\left(\bar\vartheta_\varepsilon-\vartheta _0\right)+{\rm
  I}\left(\bar\vartheta_\varepsilon\right)^{-1} \int_{0}^{T} \dot
S\left(\bar\vartheta_\varepsilon,t\right){\rm
    d}W_t\\
&\qquad + \int_{0}^{T} \frac{\dot
S\left(\bar\vartheta_\varepsilon,t\right)\left[S\left(\vartheta_0,t\right)-S\left(\bar\vartheta_\varepsilon,t\right)\right]}{
  \varepsilon{\rm
  I}\left(\bar\vartheta_\varepsilon\right) }{\rm
  d}t\\
& = \int_{0}^{T} \frac{\dot
S\left(\vartheta_,t\right)}{{\rm
  I}\left(\vartheta_0\right) }{\rm
    d}W_t -\frac{\left(\bar\vartheta_\varepsilon-\vartheta_0\right)^2 }{\varepsilon }\int_{0}^{T} \frac{\dot
S\left(\vartheta_0,t\right)\ddot S\left(\vartheta_0,t\right)}{ 2{\rm
  I}\left(\vartheta_0\right)  }{\rm
  d}t+o\left(1\right).
\end{align*} 
We see that if the seond derivative is bounded and $\varepsilon ^{-1}\varphi
_\varepsilon ^2\rightarrow 0$, then the One-step MLE is asymptotically normal
and asymptotically efficient. We say {\it formally} because there is the
problem of definition of stochastic integral. We remind this example to show
that in One-step MLE construction we need second derivative. The stochastic
models considered in the present work have no even the first derivatives and
we have to avoid the similar problem with stochastic integral too.

 Statistical problems of parameter estimation related to different
 generalizations of the model \eqref{2} with the asymptotic $\sigma =\varepsilon
 \rightarrow 0$ were treated in the works \cite{A86}, \cite{Kut88}, \cite{Kut94}
  (see review in \cite{Kut05}).  Some statistical
 problems for the related SDE with delay in the case of asymptotic
 $T\rightarrow \infty $ were studied in the works \cite{GuK99}, 
 \cite{KuK00}, \cite{Kut05}.

\section{Shift type-delay}

As it was mentioned above the problem of delay estimation we solve in two
steps. First we construct a preliminary estimator, which is consistent and
asymptotically normal and then we use this estimator and Fisher-score device
to construct One-step MLE, which is asymptotically equivalent to the MLE. In
the next section we propose and study such preliminary estimator. 

\bigskip

{\bf MDE.}
  We have observations $X^T=\left(X_t,0\leq t\leq T\right)$ of the solution of
  equation \eqref{1} and we have to estimate $\vartheta\in\Theta =\left(\alpha
  ,\beta \right)$, $0<\alpha <\beta <T$. The true value we denote as
  $\vartheta _0$. The solution $x_t$ of the equation \eqref{4} for
  $t>\vartheta $ is a function of $\vartheta $ and we write it as
  $x_t\left(\vartheta \right)$.

The MDE $\bar\vartheta _\varepsilon $ we define as solution of the following
equation
\begin{align}
\label{2.1}
\left\|X-x\left(\bar\vartheta _\varepsilon\right)\right\|=\inf_{\vartheta
  \in\Theta }\left\|X-x\left(\vartheta \right)\right\|.
\end{align}
Here  $\left\|\cdot \right\|$ is   $L_2\left[\alpha ,T\right]$-norm
\begin{align*}
\left\|X-x\left(\vartheta
\right)\right\|=\left(\int_{\alpha }^{T}\left[X_t-x_t\left(\vartheta
  \right)\right]^2{\rm d}t\right)^{1/2} .
\end{align*}
Introduce Gaussian process
$x_t^{\left(1\right)}=x_t^{\left(1\right)}\left(\vartheta _0\right) $ as
solution of the equation 
\begin{align}
\label{2.2}
{\rm d}x_t^{\left(1\right)}=S'\left(x_{t-\vartheta _0}\right)x_{t-\vartheta
  _0}^{\left(1\right)}{\rm d}t+{\rm 
  d}W_t,\qquad x_0=0, \quad 0\leq t\leq T.
\end{align}
Note that $x_t^{\left(1\right)}=\partial X_t/\partial \varepsilon
|_{\varepsilon =0} $ (see \cite{Kut94}) and that  $x_{t-\vartheta
  _0}^{\left(1\right)}=0, 0\leq t\leq \vartheta _0$. Hence
$x_t^{\left(1\right)}=W_t, 0\leq t\leq \vartheta _0$.  The derivative $\dot x_t\left(\vartheta
\right),t\geq \vartheta  $  can be calculated as follows
\begin{align}
\label{2.3}
\dot x_t\left(\vartheta \right)=\frac{\partial x_t\left(\vartheta
  \right)}{\partial \vartheta } =-\int_{\vartheta }^{t}S'\left(x_{s-\vartheta
}\right)S\left(x_{s-2\vartheta }\right){\rm d}s  .
\end{align}
Introduce as well the Gaussian random variable
\begin{align*}
\zeta \left(\vartheta _0\right)=\left(\int_{\vartheta _0}^{T}\dot
x_t\left(\vartheta_0 \right)^2{\rm d}t\right)^{-1}\int_{\vartheta
  _0}^{T}x_t^{\left(1\right)}\left(\vartheta _0\right)\dot
x_t\left(\vartheta_0 \right){\rm d}t \quad \sim\quad {\cal N}\left(0,{\rm
  D}\left(\vartheta _0\right)\right). 
\end{align*}
The variance ${\rm D}\left(\vartheta _0\right) $ of it can be calculated with
the help of the equation \eqref{2.2}, but its particular value is not
important .

\begin{theorem}
\label{T1}
Suppose that the function $S\left(x \right),x\in {\cal R}$ is positive, has
two continuous bounded derivatives and there exists $x_*\in
\left(x_0,x_{T-\beta }\right)$ such that $S'\left(x_*\right)\not=0$. Then the
MDE is consistent, asymptotically normal
\begin{align}
\label{2.4}
\varepsilon ^{-1}\left(\bar\vartheta _\varepsilon -\vartheta
_0\right)\Longrightarrow {\cal N}\left(0,{\rm D}\left(\vartheta _0\right)\right).
\end{align}
and for any $p>0$ the moments converge
\begin{align}
\label{2.5}
\varepsilon ^{-p}\Ex_{\vartheta _0} \left|\bar\vartheta _\varepsilon
-\vartheta _0\right|^p\longrightarrow  \Ex_{\vartheta _0}\left|\zeta
\left(\vartheta _0\right)\right|^p. 
\end{align}
\end{theorem}
{\bf Proof.}
First we verify the consistency of MDE.
\begin{lemma}
\label{L1}
For any $\nu >0$
\begin{align}
\label{2.6}
\sup_{\vartheta _0\in\Theta }\Pb_{\vartheta _0}\left(\left|\bar\vartheta
_\varepsilon-\vartheta _0\right|>\nu \right) \longrightarrow 0.
\end{align}
\end{lemma}
{\bf Proof.}
The proof follows standard arguments used in such problems
\begin{align*}
&\Pb_{\vartheta _0}\left(\left|\bar\vartheta
_\varepsilon-\vartheta _0\right|>\nu \right)=\Pb_{\vartheta
  _0}\left(\inf_{\left|\vartheta -\vartheta _0\right|<\nu }\left\|X-x\left(\vartheta
\right)\right\|>\inf_{\left|\vartheta -\vartheta _0\right|\geq \nu }\left\|X-x\left(\vartheta
\right)\right\|   \right)\\
&\quad \leq \Pb_{\vartheta
  _0}\left(\inf_{\left|\vartheta -\vartheta _0\right|<\nu }\left(\left\|X-x\left(\vartheta_0
\right)\right\|+\left\|x\left(\vartheta _0\right)-x\left(\vartheta
\right)\right\|\right)\right.\\
&\qquad\qquad \qquad \qquad \qquad \qquad  \left. >\inf_{\left|\vartheta -\vartheta _0\right|\geq \nu }\left(\left\|x\left(\vartheta_0
\right)-x\left(\vartheta
\right)\right\| -\left\|X-x\left(\vartheta_0
\right)\right\|   \right) \right)\\
&\quad =\Pb_{\vartheta
  _0}\left(2\left\|X-x\left(\vartheta_0
\right)\right\|>\inf_{\left|\vartheta -\vartheta _0\right|\geq \nu }\left\|x\left(\vartheta_0
\right)-x\left(\vartheta
\right)\right\| \right).
\end{align*}
Here we used the properties of norm $\left\|a+b\right\|\leq
\left\|a\right\|+\left\|b\right\|$, $\left\|a+b\right\|\geq
\left\|a\right\|-\left\|b\right\|$ and obvious equality $\inf_{\left|\vartheta
  -\vartheta _0\right|<\nu }\left\|x\left(\vartheta_0 \right)-x\left(\vartheta
\right)\right\| =0 $.

Introduce the function
\begin{align*}
g\left(\vartheta _0,\nu \right)=\inf_{\left|\vartheta -\vartheta _0\right|\geq
  \nu }\left\|x\left(\vartheta_0 \right)-x\left(\vartheta \right)\right\|
\end{align*}
and show that for any $\nu >0$ the function $g\left(\vartheta _0,\nu\right)$
satisfies $\inf_{\vartheta _0\in\Theta }g\left(\vartheta
_0,\nu\right)>0$. Suppose that $g\left(\vartheta _0,\nu\right)=0$, then there
exist $\vartheta _1$ and $\vartheta _0$ such that $\left|\vartheta
_1-\vartheta _0\right|\geq \nu $ and $\left\|x\left(\vartheta_0
\right)-x\left(\vartheta_1 \right)\right\|=0 $.  The functions
$x_t\left(\vartheta _0\right),x_t\left(\vartheta _1\right),0\leq t\leq T$ are
continuous therefore we obtain equality $x_t\left(\vartheta
_0\right)=x_t\left(\vartheta _1\right),0\leq t\leq T $. Hence for $\vartheta
_1>\theta _0$ we have equality
\begin{align*}
x_t\left(\vartheta _1\right)-x_t\left(\vartheta _0\right)=\left(\vartheta
_1-\vartheta _0\right)\dot x_t(\tilde \vartheta )\equiv 0,\qquad 0\leq t\leq T,
\end{align*}
where $\vartheta _0\leq \tilde \vartheta \leq \vartheta _1$. As the function
$S\left(x\right)$ is strictly positive, this equality implies (see
\eqref{2.3}) $S'\left(x_{s-\tilde \vartheta }\right)\equiv 0 $ for all $s\in
\left[\vartheta _0\leq s\leq T\right]$. It is possible to take such $s $ that
$s-\tilde\vartheta =x_*$, but by condition of Theorem $S'\left(x_{s-\tilde
  \vartheta }\right)=S'\left(x_* \right)\not= 0 $. Hence $g\left(\vartheta
_0,\nu \right)>0$.

By Tchebychev inequality
\begin{align*}
\Pb_{\vartheta _0}\left(\left|\bar\vartheta _\varepsilon-\vartheta
_0\right|>\nu \right)\leq \frac{4}{g\left(\vartheta _0,\nu
\right)^2 }\Ex_{\vartheta _0}\int_{\alpha }^{T}\left[X_t-x_t\left(\vartheta
  _0\right)\right]^2{\rm d}t \leq  \frac{C\;\varepsilon ^2}{g\left(\vartheta _0,\nu
\right)^2 }\rightarrow 0.
\end{align*}
Here we used the estimate \eqref{A2} proved in Appendix below. 
Therefore \eqref{2.5} is proved.

To prove the asymptotic normality \eqref{2.4} we note that the estimator
$\bar\vartheta _\varepsilon $ satisfies the minimum distance equation (MDEq)
\begin{align*}
\int_{\alpha }^{T}\left[X_t-x_t\left(\bar\vartheta _\varepsilon
  \right)\right]\dot x_t\left(\bar\vartheta _\varepsilon
  \right){\rm d}t=0.
\end{align*}
Let us introduce the random variable  $\bar u_\varepsilon =\varepsilon
^{-1}\left(\bar\vartheta 
_\varepsilon-\vartheta _0 \right)$  and write the expansion  
\begin{align*}
x_t\left(\bar\vartheta _\varepsilon \right)&=x_t\left(\vartheta _0+\varepsilon 
\bar u _\varepsilon \right)=x_t\left(\vartheta _0\right)+\varepsilon 
\bar u _\varepsilon\dot  x_t\left(\vartheta
_0\right)+O\left(\varepsilon ^2\right),\\
\dot x_t\left(\bar\vartheta _\varepsilon \right)&=\dot x_t\left(\vartheta
_0\right)   +O\left(\varepsilon \right).
\end{align*}
 Then the MDEq allows us to write
\begin{align*}
\int_{\alpha }^{T}\varepsilon ^{-1} \left( X_t-x_t\left(\vartheta
  _0\right) \right)\left[\dot
x_t\left(\vartheta _0 \right)+O\left(\varepsilon
\right)\right]{\rm d}t=\bar u 
_\varepsilon\int_{\alpha }^{T} \dot 
x_t\left(\vartheta _0 \right)^2{\rm d}t \left(1+O\left(\varepsilon \right)\right)
\end{align*}
and
\begin{align*}
\bar u _\varepsilon=\left(\int_{\alpha }^{T} \dot 
x_t\left(\vartheta _0 \right)^2{\rm d}t \right)^{-1}\int_{\alpha }^{T}
\varepsilon ^{-1} \left( X_t-x_t\left(\vartheta 
  _0\right) \right)\dot 
x_t\left(\vartheta _0 \right){\rm d}t.
\end{align*}
Here
\begin{align*}
\varepsilon ^{-1} \left( X_t-x_t\left(\vartheta 
  _0\right) \right)\longrightarrow x_t^{\left(1\right)}\left(\vartheta _0\right)
\end{align*}
in probability and therefore
\begin{align*}
\bar u _\varepsilon=\varepsilon
^{-1}\left(\bar\vartheta 
_\varepsilon-\vartheta _0 \right)\Longrightarrow \zeta \left(\vartheta _0\right).
\end{align*}
Moreover, it is possible to show that we have asymptotic normality in
probability, because the limit random variable $  \zeta \left(\vartheta
_0\right)$ is defined on the same probability space. 

For $\left|\vartheta -\vartheta
_0\right|\leq \nu $ and  sufficiently small $\nu >0$   we have
\begin{align*}
\left\|x\left(\vartheta \right)-x\left(\vartheta_0
\right)\right\|^2&=\left(\vartheta-\vartheta_0\right)^2 \|\dot
x(\tilde\vartheta )\|^2 =\left(\vartheta-\vartheta_0\right)^2 \|\dot
x(\vartheta_0 )\|^2\left(1+o\left(1\right)\right)\\
&\geq
\frac{1}{2}\left(\vartheta-\vartheta_0\right)^2 \|\dot x(\vartheta_0 )\|^2.
\end{align*}
This estimate together with $g\left(\vartheta _0,\nu \right)>0$ allow us to
write the estimate: there exists $\kappa >0$ such that 
\begin{align*}
g\left(\vartheta _0,\nu \right)\geq \kappa \nu .
\end{align*}
Let us put $\nu =\varepsilon y>0$, then for any integer $N>1$
\begin{align*}
&\Pb_{\vartheta _0}\left(\varepsilon ^{-1} \left|\bar\vartheta _\varepsilon
-\vartheta _0\right|>y\right) \leq \frac{2^{2N}}{\kappa
  ^{2N}y^{2N}}\Ex_{\vartheta _0}\left|\int_{\alpha
}^{T}\left|X_t-x_t\left(\vartheta _0\right)\right|^2{\rm d}t\right|^N\\
&\qquad \qquad \qquad \leq \frac{2^{2N}\left(T-\alpha \right)^{N-1}}{\varepsilon ^{2N}\kappa
  ^{2N}y^{2N}}\int_{\alpha
}^{T}\Ex_{\vartheta _0}\left|X_t-x_t\left(\vartheta _0\right)\right|^{2N}{\rm
  d}t\leq \frac{C_N}{y^{2N}}, 
\end{align*}
where the constant $C>0$ does not depend on $\varepsilon $ and $\vartheta _0$.

Let us denote $F_\varepsilon \left(\vartheta _0,y\right)=\Pb_{\vartheta _0}\left(\left|\bar
u_\varepsilon \right|<y\right)$. Then for any $p>0$ we take $2N\geq p+1$ and can
write
\begin{align*}
\Ex_{\vartheta _0}\left|\frac{\bar\vartheta _\varepsilon -\vartheta _0
}{\varepsilon }\right|^p&=\int_{0}^{\infty }y^p{\rm d}
F_\varepsilon\left(\vartheta _0,y\right)=1- \int_{1}^{\infty }y^p{\rm
  d}\left(1-F_\varepsilon\left(\vartheta _0,y\right)\right)\\ 
&\leq
2+p\int_{1}^{\infty }y^{p-1}\left(1-F_\varepsilon\left(\vartheta
_0,y\right)\right){\rm d}y \\
&=
2+p\int_{1}^{\infty }y^{p-1}\Pb_{\vartheta _0}\left(\varepsilon ^{-1}
\left|\bar\vartheta _\varepsilon -\vartheta _0\right|>y\right){\rm d}y \\
&\leq 2+pC_p\int_{1}^{\infty }y^{p-1-2N}{\rm d}y\leq C.
\end{align*}
Therefore for any $p>1$ the random variables $\varepsilon ^{-p}\left|\bar\vartheta
_\varepsilon -\vartheta _0\right|^p$ are uniformly integrable 
and we obtain the convergence of moments \eqref{2.5}.

Remark that if the value $x_*$ mentioned in the Theorem \ref{T1} does not
exist and $S'\left(x\right)=0$ for all $x\in \left[x_0,x_{T-\beta }\right] $,
then the function $S\left(x\right), x\in \left[x_0,x_{T-\beta }\right]$ does
not depend on $x$ and consistent estimation of $\vartheta _0$ is impossible. 

\bigskip

{\bf One-step MLE.}
Let us denote $\left\{\Pb_{\vartheta }^{\left(\varepsilon \right)},\vartheta
\in \Theta \right\}$ the family of measures induced on the measurable space
$\left({\cal C}_{\left(0,T\right)},{\scr B}\right)$ of continuous on $\left[0,T\right] $
functions by the solutions of the equation \eqref{1} with different
$\vartheta \in\Theta $. These measures are equivalent with the likelihood
ratio function $L\left(\vartheta ,X^T\right)$ \cite{LS01}. Introduce the
normalized likelihood ratio process
\begin{align*}
Z_\varepsilon \left(u\right)=\frac{L\left(\vartheta_0+\varepsilon u
  ,X^T\right)}{L\left(\vartheta_0 ,X^T\right)},\qquad u\in \UU_\varepsilon
=\left(\frac{\alpha -\vartheta _0}{\varepsilon },\frac{\beta  -\vartheta
  _0}{\varepsilon }\right). 
\end{align*}
We have (see \eqref{LR})
\begin{align*}
\ln Z_\varepsilon \left(u\right)&=\int_{0}^{T}\frac{S\left(X_{t-\vartheta
    _0-\varepsilon u}\right)-S\left(X_{t-\vartheta
    _0}\right)}{\varepsilon }\;{\rm d}W_t\\
&\qquad - \int_{0}^{T}\frac{\left[S\left(X_{t-\vartheta
    _0-\varepsilon u}\right)-S\left(X_{t-\vartheta
    _0}\right)\right]^2}{2\varepsilon^2 }\;{\rm d}t.
\end{align*}
As it follows from the given below proof of the theorem (see also \cite{A86}),
this likelihood ratio admits the representation
\begin{align*}
Z_\varepsilon \left(u\right)&=\exp\left\{-u\int_{\vartheta _0}^{T}S'\left(x_{t-\vartheta
    _0}\right)S\left(x_{t-2\vartheta
    _0}\right)\;{\rm d}W_t-\frac{u^2}{2}{\rm
  I}\left(\vartheta_0 \right)+o\left(1\right) \right\}.
\end{align*}
Here $x_{t-\vartheta
    _0}=x_{t-\vartheta
    _0}\left(\vartheta _0\right) $, $x_{t-2\vartheta
    _0}=x_{t-2\vartheta
    _0}\left(\vartheta _0\right) $ and ${\rm
  I}\left(\vartheta \right) $ plays the role of   Fisher information.  It is
given by the equality
\begin{align*}
 {\rm I}\left(\vartheta \right)=\int_{\vartheta }^{T}S'\left(x_{t-\vartheta
 }\left(\vartheta \right)\right)^2S\left(x_{t-2\vartheta
 }\left(\vartheta\right)\right)^2 {\rm d}t.
\end{align*}
As
\begin{align*}
\int_{\vartheta_0 }^{T}S'\left(x_{t-\vartheta
    _0}\right)S\left(x_{t-2\vartheta
    _0}\right)\;{\rm d}W_t\quad \sim\quad {\cal N}\left(0,{\rm
  I}\left(\vartheta_0 \right)\right)
\end{align*}
the family of measures $\left\{\Pb_{\vartheta }^{\left(\varepsilon \right)},\vartheta
\in \Theta \right\} $ is locally asymptotically normal (LAN) for all
$\vartheta _0\in\Theta $. Therefore we have Hajek-Le Cam's lower bound on the
mean squared 
risk of all estimators $\vartheta _\varepsilon $
\begin{align*}
\lim_{\nu \rightarrow 0}\Liminf_{\varepsilon \rightarrow
  0}\sup_{\left|\vartheta -\vartheta _0\right|\leq \nu }\varepsilon
^{-2}\Ex_\vartheta \left(\vartheta _\varepsilon -\vartheta \right)^2\geq {\rm 
  I}\left(\vartheta_0 \right)^{-1}.
\end{align*}
As usual, we call estimator $\vartheta _\varepsilon ^*$ asymptotically
efficient if for this estimator we have equality
\begin{align}
\label{HLC}
\lim_{\nu \rightarrow 0}\lim_{\varepsilon \rightarrow 0}\sup_{\left|\vartheta
  -\vartheta _0\right|\leq \nu }\varepsilon ^{-2}\Ex_\vartheta \left(\vartheta
_\varepsilon^* -\vartheta \right)^2= {\rm I}\left(\vartheta_0 \right)^{-1}
\end{align}
for all $\vartheta _0\in\Theta $.
 Recall that the MLE $\hat\vartheta _\varepsilon $ for this model of
 observations is asymptotically efficient \cite{A86}, \cite{Kut94}.

One-step MLE $\vartheta _\varepsilon ^\star $ for this model of observations
formally can be written as follows
\begin{align*}
\vartheta _\varepsilon ^\star&=\bar\vartheta _\varepsilon +\varepsilon ^2{\rm
  I}\left(\bar\vartheta _\varepsilon\right)^{-1} \left.\frac{\partial \ln
  L\left(\vartheta ,X^T\right)}{\partial \vartheta }\right|_{\vartheta
  =\bar\vartheta _\varepsilon}\\
&=\bar\vartheta _\varepsilon +{\rm
  I}\left(\bar\vartheta _\varepsilon\right)^{-1}
\int_{0}^{T}S'\left(X_{t-\bar\vartheta _\varepsilon}\right) \frac{\partial
  X_{t-\bar\vartheta _\varepsilon}}{\partial \vartheta }\left[{\rm
    d}X_t-S\left(X_{t-\bar\vartheta _\varepsilon}\right){\rm d}t \right] .
\end{align*}

We say {\it formally} because in this writing there are two difficulties. The
first one: as the estimator $\bar\vartheta _\varepsilon $ depends on all
observations on the interval $\left[\alpha ,T\right]$ the It\^o stochastic
integral is not well defined.  The second problem, of course, is the
calculation of the derivative ${\partial
  X_{t-\bar\vartheta _\varepsilon}}/{\partial \vartheta } $. Our goal is to
find an alternative expression for One-step MLE which does not have such
problems. 

Let us write the expression for $X_{t-\vartheta }, $ without Wiener part
\begin{align*}
\hat X_{t-\vartheta }=x_0+\int_{0}^{t-\vartheta }S\left(X_{s-\vartheta
  _0}\right)\,{\rm d}s .
\end{align*}
Then 
\begin{align*}
\frac{\partial \hat 
  X_{t-\vartheta}}{\partial \vartheta }=-S\left(X_{t-\vartheta -\vartheta
  _0}\right)\1_{\left\{ t> \vartheta \right\}}.
\end{align*}

First we replace the stochastic process $X_{t-\bar\vartheta _\varepsilon }$ by
$x_{t-\bar\vartheta _\varepsilon }\left(\bar\vartheta _\varepsilon \right)$
and write the integral as follows
\begin{align*}
&\int_{0}^{T}S'\left(x_{t-\bar\vartheta _\varepsilon}\right) \left.\frac{\partial
    x_{t-\vartheta }\left(\bar\vartheta _\varepsilon
    \right)}{\partial \vartheta }\right|_{\vartheta =\bar\vartheta _\varepsilon }\left[{\rm
      d}X_t-S\left(X_{t-\bar\vartheta _\varepsilon}\right){\rm d}t
    \right]\\ 
&\quad = \int_{\bar\vartheta _\varepsilon }^{T}S'\left(x_{t-\bar\vartheta _\varepsilon}\right)
  S\left(x_{t-2\bar\vartheta _\varepsilon}\right) S\left(X_{t-\bar\vartheta
    _\varepsilon}\right){\rm d}t -\int_{\bar\vartheta _\varepsilon}^{T}S'\left(x_{t-\bar\vartheta
    _\varepsilon}\right)S\left(x_{t-2\bar\vartheta _\varepsilon}\right) {\rm
    d}X_t.
\end{align*}
Note that $X_{t-\vartheta }|_{\varepsilon =0}=x_{t-\vartheta }\left(\vartheta
_0\right)$ and we differentiate on $\vartheta $ only, but then replace the
second value $\vartheta_0 $ by $\bar\vartheta _\varepsilon $ too.  Remind that
the likelihood ratio depends on the true value $\vartheta _0$ but we never
differentiate it on $\vartheta _0$.

 Further, let us denote
\begin{align*}
H\left(t,\vartheta \right)=S'\left(x_{t-\vartheta }\left(\vartheta \right)\right)
S\left(x_{t-2\vartheta }\left(\vartheta \right)\right).  
\end{align*}
Then
\begin{align*}
{\rm d}\left[H\left(t,\vartheta \right)X_t \right]=H\left(t,\vartheta
\right){\rm d}X_t  +X_tH'_t\left(t,\vartheta \right){\rm d}t
\end{align*}
and
\begin{align*}
\int_{\vartheta  }^{T}S'\left(x_{t-\vartheta }\right) S\left(x_{t-2\vartheta
}\right){\rm d}X_t=H\left(T,\vartheta \right)X_T-H\left(\vartheta ,\vartheta
\right)X_\vartheta -\int_{\vartheta }^{T} H'_t\left(t,\vartheta \right)X_t {\rm d}t.
\end{align*}
Of course, we have
\begin{align*}
H'_t\left(t,\vartheta \right)&=S''\left(x_{t-\vartheta }\left(\vartheta
\right)\right) S\left(x_{t-2\vartheta }\left(\vartheta
\right)\right)^2\\
&\qquad \qquad \qquad +S'\left(x_{t-\vartheta }\left(\vartheta \right)\right)
S'\left(x_{t-2\vartheta }\left(\vartheta \right)\right)S\left(x_{t-3\vartheta
}\left(\vartheta \right)\right)\1_{\left\{t\geq 2\vartheta \right\}}.
\end{align*}
 Now we can put
\begin{align*}
\Psi \left(\bar\vartheta _\varepsilon \right)=H\left(T,\bar\vartheta
_\varepsilon \right)X_T-H\left(\bar\vartheta _\varepsilon ,\bar\vartheta _\varepsilon
\right)X_{\bar\vartheta _\varepsilon} -\int_{\bar\vartheta _\varepsilon }^{T}
H'_t\left(t,\bar\vartheta _\varepsilon \right)X_t\; {\rm d}t 
\end{align*}
and define the One-step MLE
\begin{align}
\label{2.7}
\vartheta _\varepsilon ^\star&=\bar\vartheta _\varepsilon +{\rm I}
\left(\bar\vartheta _\varepsilon\right)^{-1} \left[
  \int_{\bar\vartheta _\varepsilon }^{T}S'\left(x_{t-\bar\vartheta _\varepsilon}\right)
  S\left(x_{t-2\bar\vartheta _\varepsilon}\right) S\left(X_{t-\bar\vartheta
    _\varepsilon}\right){\rm d}t -\Psi \left(\bar\vartheta _\varepsilon
  \right)\right] .
\end{align}
Its properties are given in the next theorem.
\begin{theorem}
\label{T2}
Suppose that the function $S\left(\cdot \right)$ is positive,  has four
continuous bounded derivatives  and  there exists $x_*\in \left(x_0,x_{T-\beta
}\right)$ such that  $S'\left(x_*\right)\not=0$. Then the
One-step MLE $\vartheta _\varepsilon ^\star$ is consistent, asymptotically
normal
\begin{align}
\label{2.8}
\varepsilon ^{-1}\left(\vartheta _\varepsilon ^\star-\vartheta
_0\right)\Longrightarrow  {\cal N}\left(0,{\rm I}
\left(\vartheta _0\right)^{-1}\right).
\end{align}
 and asymptotically efficient. 
\end{theorem}
{\bf Proof.} Let us study the expression
\begin{align*}
R\left(\bar\vartheta _\varepsilon \right)&=\int_{\bar\vartheta _\varepsilon }^{T}S'\left(x_{t-\bar\vartheta _\varepsilon}\right)
  S\left(x_{t-2\bar\vartheta _\varepsilon}\right) S\left(X_{t-\bar\vartheta
    _\varepsilon}\right){\rm d}t -\Psi \left(\bar\vartheta _\varepsilon
  \right)= H\left(\bar\vartheta _\varepsilon ,\bar\vartheta _\varepsilon
\right)X_{\bar\vartheta _\varepsilon} \\
&\quad + \int_{\bar\vartheta _\varepsilon }^{T}\left[H\left(t,\bar\vartheta
_\varepsilon \right) S\left(X_{t-\bar\vartheta
    _\varepsilon}\right)-H'_t\left(t,\bar\vartheta _\varepsilon \right)X_t\right]{\rm d}t-   H\left(T,\bar\vartheta
_\varepsilon \right)X_T.
\end{align*}
 We have to expand it at the vicinity of the point $\vartheta _0$.  As $\bar\vartheta
_\varepsilon =\vartheta _0+\varepsilon \bar u_\varepsilon $ and the function
$H\left(t,\vartheta \right)$ is sufficiently smooth w.r.t. $\vartheta $ and
$t$ we can write the expansions for the function $H\left(\cdot \right)$.
 Before we note that  there is a problem with expansion of the
term  $H(\bar \vartheta _\varepsilon ,\bar \vartheta _\varepsilon )X_{\bar
  \vartheta _\varepsilon }$. The process $X_t$ has no derivative but we have 
\begin{align*}
X_{\bar \vartheta _\varepsilon }&=X_{\vartheta _0}+\int_{\vartheta _0}^{\bar
  \vartheta _\varepsilon } S\left(X_{t-\vartheta _0}\right){\rm d}t
+\varepsilon \left[ W_{\bar \vartheta _\varepsilon }-W_{\vartheta
    _0}\right]\\
&=X_{\vartheta _0}+\varepsilon \bar u_\varepsilon S\left(x_0\right)
\left(1+o\left(1\right)\right)+\varepsilon 
^{3/2}O\left(1\right) .
\end{align*}
Therefore we can write
\begin{align*}
H\left(T,\bar\vartheta _\varepsilon \right)&=H\left(T,\vartheta _0
\right)+\varepsilon \bar u_\varepsilon H'_\vartheta
(T,\vartheta _0 )+ \varepsilon^2  A
(T,\tilde\vartheta _\varepsilon ),\\
H\left(\bar\vartheta_\varepsilon ,\bar\vartheta_\varepsilon \right)
&=H\left(\vartheta _0,\vartheta _0 \right)+\varepsilon \bar
u_\varepsilon \left[H'_t (\vartheta _0 ,\vartheta
_0)+H'_\vartheta (\vartheta _0 ,\vartheta
_0)\right]+ {\varepsilon^2 } B
(\tilde \vartheta _\varepsilon  ,\tilde\vartheta _\varepsilon ),
\end{align*}
and for the integral we have 
\begin{align*}
 &\int_{\bar\vartheta_\varepsilon }^{T} H'_t\left(t,\bar\vartheta
_\varepsilon \right)X_t {\rm d}t=\int_{\vartheta _0 }^{T}
H'_t\left(t,\vartheta _0 \right)X_t {\rm d}t\\ 
&\quad  +\varepsilon \bar
u_\varepsilon \left[\int_{\bar\vartheta_\varepsilon}^{T} H''_{t,\vartheta
}(t,\vartheta _0)X_t {\rm d}t-H'_t\left(\vartheta _0,\vartheta _0
\right)X_{\vartheta _0}\right]  +\varepsilon^{3/2}C(T,\tilde \vartheta _\varepsilon  ,\tilde\vartheta
_\varepsilon  ). 
\end{align*}
Here we denoted $A\left(\cdot \right),B\left(\cdot \right)$ and $C\left(\cdot
\right)$ the corresponding residuals. Note that all polynomial moments of these
quantities are finite.  Of course, the values of $\tilde
\vartheta _\varepsilon $ in these three expansions are different and this is
just a symbolical writing. 

Therefore we have the presentation
\begin{align*}
\Psi \left(\bar \vartheta _\varepsilon \right)=\Psi \left( \vartheta _0
\right)+\varepsilon \bar
u_\varepsilon \hat \Psi
(\vartheta _0 ) + \varepsilon^{3/2} O (1 )   ,
\end{align*}
where 
\begin{align*}
\Psi \left( \vartheta _0 \right)&=\int_{\vartheta_0
}^{T}S'\left(x_{t-\vartheta_0 }\right) S\left(x_{t-2\vartheta_0 }\right){\rm
  d}X_t\\
 &=\int_{\vartheta_0 }^{T}H\left(t,\vartheta
_0\right)S\left(X_{t-\vartheta_0 }\right){\rm d}t+\varepsilon
\int_{\vartheta_0 }^{T}H\left(t,\vartheta _0\right){\rm d}W_t
\end{align*}
and 
\begin{align*}
&\hat \Psi (\vartheta _0 )= H'_\vartheta \left(T,\vartheta _0\right)X_T-
  \left[ H'_\vartheta\left(\vartheta _0,\vartheta _0\right)+
    H_t'\left(\vartheta _0,\vartheta _0\right)\right]X_{\vartheta _0}-
  H\left(\vartheta _0,\vartheta _0\right)S\left(x_0\right) \\
 &\qquad 
  \qquad+H_t'\left(\vartheta _0,\vartheta _0\right)X_{\vartheta _0}
  -\int_{\vartheta _0}^{T} H_{t,\vartheta }''\left(t,\vartheta
  _0\right)X_{t}{\rm d}t \\ 
&\qquad  = H'_\vartheta \left(T,\vartheta
  _0\right)X_T- H'_\vartheta\left(\vartheta _0,\vartheta _0\right)X_{\vartheta
    _0}  - H\left(\vartheta _0,\vartheta
  _0\right)S\left(x_0\right)\\
&\qquad \qquad -\int_{\vartheta _0}^{T} H''_{t,\vartheta }\left(t,\vartheta
  _0\right)X_{t}{\rm d}t\\
 &\qquad  =\int_{\vartheta _0}^{T} H'_\vartheta
  \left(t,\vartheta _0\right){\rm d}X_{t} - H\left(\vartheta _0,\vartheta
  _0\right)S\left(x_0\right)\\
 &\qquad =\int_{\vartheta _0}^{T} H'_\vartheta
  \left(t,\vartheta _0\right)S\left(X_{t-\vartheta _0}\right){\rm
    d}t+\varepsilon \int_{\vartheta _0}^{T} H'_\vartheta\left(t,\vartheta
  _0\right){\rm d}W_{t} - H\left(\vartheta _0,\vartheta
  _0\right)S\left(x_0\right).
\end{align*}
Further
\begin{align*}
\int_{\bar\vartheta _\varepsilon }^{T} H\left(t,\bar\vartheta
_\varepsilon\right)S\left(X_{t-\bar\vartheta _\varepsilon}\right){\rm d}t&
=\int_{\vartheta _0}^{T} H\left(t,\vartheta _0\right)S\left(X_{t-\bar\vartheta
  _\varepsilon}\right){\rm d}t -\varepsilon \bar u_\varepsilon
H\left(\vartheta _0,\vartheta _0\right)S\left(x_0\right)\\
 &\qquad
+\varepsilon \bar u_\varepsilon \int_{\vartheta _0}^{T} H'_\vartheta
\left(t,\vartheta _0\right)S\left(X_{t-\vartheta _0}\right){\rm
  d}t+O\left(\varepsilon ^{3/2}\right).
\end{align*}
Therefore
\begin{align*}
R\left(\bar\vartheta _\varepsilon \right)&=\int_{\vartheta _0}^{T} H\left(t,\vartheta _0\right)\left[S\left(X_{t-\bar\vartheta
  _\varepsilon}\right)-S\left(X_{t-\vartheta
  _0}\right) \right]{\rm d}t\\
&\qquad \qquad -\varepsilon \int_{\vartheta _0}^{T}
H\left(t,\vartheta _0\right){\rm d}W_t +O\left(\varepsilon ^{3/2}\right)\\
&=-\varepsilon \bar u_\varepsilon \int_{\vartheta _0}^{T}
 H\left(t,\vartheta _0\right)S'\left(x_{t-\vartheta
  _0}\right)S\left(x_{t-\vartheta
  _0}\right){\rm d}t\\
&\qquad \qquad-\varepsilon \int_{\vartheta _0}^{T}
H\left(t,\vartheta _0\right){\rm d}W_t +O\left(\varepsilon ^{3/2}\right)\\
&=-\varepsilon \bar u_\varepsilon{\rm I}\left(\vartheta _0\right)-\varepsilon
\int_{\vartheta _0}^{T} 
H\left(t,\vartheta _0\right){\rm d}W_t +O\left(\varepsilon ^{3/2}\right).
\end{align*}
The obtained relations allow us to write 
\begin{align*}
\frac{\vartheta _\varepsilon ^\star-\vartheta _0}{\varepsilon }&= \bar
u_\varepsilon-{\rm I}
\left(\bar\vartheta _\varepsilon\right)^{-1} \int_{\vartheta_0 }^{T}
H\left(t,\vartheta _0\right){\rm d}W_t-\bar u_\varepsilon \frac{{\rm I}\left(\vartheta _0\right)}{{\rm I}
\left(\bar\vartheta _\varepsilon\right)} +O\left(\varepsilon ^{1/2}\right)\\
&=-{\rm I}
\left(\vartheta _0\right)^{-1} \int_{\vartheta _0 }^{T}H\left(t,\vartheta
_0\right){\rm d}W_t \left(1+o\left(\varepsilon ^{1/2}\right)\right).
\end{align*}
Recall that
\begin{align*}
-{\rm I}
\left(\vartheta _0\right)^{-1} \int_{\vartheta _0 }^{T}H\left(t,\vartheta
_0\right){\rm d}W_t \quad \sim\quad {\cal N}\left(0,{\rm I}
\left(\vartheta _0\right)^{-1} \right). 
\end{align*}
Note that the Fisher information is separated from zero, i.e.,
\begin{align*}
\inf_{\vartheta \in\Theta }{\rm I}\left(\vartheta \right)>0
\end{align*}
and we can easily verify the convergence of all polynomial moments of
$o\left(\varepsilon ^{1/2}\right)$ above. Moreover, in all  estimates 
of the type
\begin{align*}
\Ex_{\vartheta _0}\left|o\left(\varepsilon ^{1/2}\right)\right|^p\leq
C\varepsilon ^{\frac{p}{2}} 
\end{align*}
with any $p>0$ and corresponding  constants
$C=C\left(K\right)>0$ not 
depending on $\vartheta _0\in K$ for any compact set  $K\subset\Theta $.  This
provides us the uniform convergence of the moments of $u_\varepsilon
^\star=\varepsilon ^{-1}\left( \vartheta _\varepsilon ^\star-\vartheta
_0\right)$. For the second moment we have 
\begin{align*}
\sup_{\left|\vartheta -\vartheta _0\right|\leq \nu }\Ex_\vartheta
\left|\frac{\vartheta _\varepsilon ^\star-\vartheta}{\varepsilon
}\right|^2\xrightarrow{\ \varepsilon \rightarrow 0\ } \sup_{\left|\vartheta
  -\vartheta _0\right|\leq \nu } {\rm I}\left(\vartheta
\right)^{-1}\xrightarrow{\ \nu \rightarrow 0\ }{\rm I}\left(\vartheta_0
\right)^{-1}.
\end{align*}
 
Therefore the One-step MLE $\vartheta _\varepsilon ^\star $ is asymptotically
efficient.

\bigskip

{\bf Example 1}. Suppose that we have Ornstein-Uhlenbeck process with
delay
\begin{align*}
{\rm d}X_t=-\gamma X_{t-\vartheta }\;{\rm d}t+\varepsilon\, {\rm d}W_t,\qquad
X_s=x_0,\quad s\leq 0,\quad 0\leq t\leq T.
\end{align*}
Here $\gamma >0$ is supposed to be known, $x_0>0$ and $\vartheta \in\left(\alpha,\beta
\right)$, $0<\alpha <\beta <T$. The solution $x_t=x_t\left(\vartheta
_0\right)$ of the limit  equation 
\begin{align*}
\frac{{\rm d}x_t}{{\rm d}t}=-\gamma x_{t-\vartheta_0 },\qquad x_s=x_0,\; s\leq
0,\qquad 0\leq t\leq T 
\end{align*}
 can be written
in explicit form as follows
\begin{align*}
x_t\left(\vartheta _0\right)=x_0y_t\left(\vartheta _0\right)-\gamma x_0 \int_{-\vartheta _0}^{0}y_{t-s-\vartheta _0}\left(\vartheta _0\right){\rm d}s,\qquad t\geq 0,
\end{align*}
where
\begin{align*}
y_t\left(\vartheta _0\right)=\sum_{k=0}^{\left[\frac{t}{\vartheta _0}\right]}\left(-1\right)^k\frac{\gamma ^k}{k!}\left(t-k\vartheta _0\right)^k,
\end{align*}
$\left[a\right]$ is entire part of $a$ (see \cite{KuM92} for details). 

The preliminary MDE $\bar\vartheta _\varepsilon $ is
\begin{align*}
\bar\vartheta _\varepsilon=\arg\inf_{\alpha <\vartheta <\beta  }\int_{\alpha
}^{T}\left[X_t-x_t\left(\vartheta \right)\right]^2{\rm d}t .
\end{align*}
 Recall that the Fisher information in this
problem is
\begin{align*}
{\rm I}\left(\vartheta \right)=\gamma ^4\int_{\vartheta }^{T}x_{t-2\vartheta
}\left(\vartheta \right)^2{\rm d}t .
\end{align*} 
The One-step MLE $\vartheta _\varepsilon ^\star$ has the following form
\begin{align*}
\vartheta _\varepsilon ^\star&=\bar\vartheta _\varepsilon +{\rm
  I}\left(\bar\vartheta _\varepsilon \right)^{-1} \Bigl[\gamma
  ^2x_0X_{\bar\vartheta _\varepsilon} 
  -\gamma ^2x_{T-2\bar\vartheta _\varepsilon }(\bar\vartheta
  _\varepsilon)X_T\\
& \qquad \qquad   -\gamma ^3\int_{\bar\vartheta _\varepsilon}^{T}
  x_{t-2\bar\vartheta _\varepsilon}(\bar\vartheta
  _\varepsilon)X_{t-\bar\vartheta
      _\varepsilon}{\rm d}t-\gamma ^3\int_{2\bar\vartheta _\varepsilon}^{T}   x_{t-3\bar\vartheta _\varepsilon}(\bar\vartheta
  _\varepsilon)X_{t}{\rm d}t\Bigr] . 
\end{align*}
By Theorem \ref{T2} we have
\begin{align*}
\varepsilon ^{-1}\left(\vartheta _\varepsilon ^\star-\vartheta
_0\right)\Longrightarrow  {\cal N}\left(0,{\rm I}\left(\vartheta \right)^{-1}
\right). 
\end{align*}

\section{Scale-type delay}

Below we consider the models of type \eqref{1a}. We start with the linear
model with the equation, which can be  called {\it Stochastic Pantograph Equation }
and then we  study the nonlinear model.

\bigskip

{\bf Stochastic Pantograph Equation.}
 Recall that the deterministic {\it Pantograph
equation}
\begin{align}
\label{14a}
\frac{{\rm d}x_t}{{\rm d}t}=a\,x_t+b\,x_{\vartheta t},\qquad x_0,\quad t\geq 0,
\end{align}
was introduced by Ockendon and Tayler in 1971 \cite{OT71}. This equation
describes the special construction on the electric locomotive of the British
Railways, which allowed to collect the current from an overhead wire.  Here
the parameter $0<\vartheta <1$ defines the delay $d=t-\vartheta t$. The
further study of this equation can be found in the thesis \cite{TG17}.

Introduce the condition

\bigskip

${\scr A}$. {\it The reals $x_0\not=0$, $b\not=0$ and $a+b\not=0$.}

\bigskip

Note that if $x_0=0$, then the solution of \eqref{14a} is $x_t\equiv 0$. If
$b=0$, then there is no delay. If $a+b=0$, then we have once more a trivial
solution $x_t\equiv x_0$.

The stochastic version of Pantograph equation can be 
\begin{align}
\label{14}
{\rm d}X_t=a\,X_t\,{\rm d}t+b\,X_{\vartheta t}\,{\rm d}t+\varepsilon \,  {\rm
  d}W_t,\qquad X_0=x_0,\qquad t\geq 0.
\end{align}
Consider the problem of estimation $\vartheta \in\Theta =\left(\alpha ,\beta
\right)$, $0<\alpha <\beta <1$ by the  observations $X^T=\left(X_t,0\leq t\leq
T\right)$. The reals $a, b\not=0,x_0\not=0$ and $\varepsilon \in (0,1]$ are
  supposed to be known. The likelihood ratio function is
\begin{align*}
L\left(\vartheta ,X^T\right)=\exp\left\{\int_{0}^{T}\frac{bX_{\vartheta t}}{\varepsilon
  ^2}{\rm d}X_t -\int_{0}^{T} \frac{\left[b^2X_{\vartheta t}^2-2abX_{
      t}X_{\vartheta t}\right]}{2\varepsilon ^2}{\rm d}t\right\},\quad
\vartheta \in \Theta 
\end{align*} 
and numerical calculation of the MLE $\hat \vartheta _\varepsilon $ defined by
the equation \eqref{mle} has the same problems as in the case of observations
\eqref{1}. Note that the asymptotic ($\varepsilon \rightarrow 0$)  properties
of the MLE in the case $a=0$ were described in \cite{Kut88}. Here we consider
the construction of the one-step MLE using the Fisher-score device similar to
that given above. 

The role of  Fisher information  plays the quantity
 \begin{align*}
{\rm I}\left(\vartheta \right)=b^2 \int_{0}^{T}t^2 \left[ax_{\vartheta   
  t}\left(\vartheta \right)+bx_{\vartheta^2
  t}\left(\vartheta \right)\right]^2 {\rm d}t.
\end{align*}
As preliminary estimator wa can take the MDE $\bar\vartheta _\varepsilon $
defined by the equation \eqref{2.1}, where $x\left(\vartheta
\right)=\left(x_t\left(\vartheta \right),0\leq t\leq T\right)$ is the solution 
of the equation \eqref{14}. It can be shown that this estimator is consistent 
and asymptotically normal $\varepsilon ^{-1}\left(\bar\vartheta
_\varepsilon-\vartheta _0\right)$. Moreover, it is consistent and 
asymptotically normal even if we take observations $X^{\tau _\varepsilon
}=\left(X_t,0\leq t\leq \tau _\varepsilon \right)$ on the interval
$\left[0,\tau _\varepsilon \right] $ with $\tau _\varepsilon \rightarrow
0$.

\bigskip

{\bf Estimator of substitution.} 
Note that the model \eqref{14} allows us to use another more  simple estimator 
as preliminary one, which we introduce as follows  
\begin{align}
\label{17a}
\bar\vartheta _{\tau_\varepsilon }
=\frac{2}{\left(a+b\right)b\tau_{\varepsilon } ^2}\left[X_{\tau _\varepsilon
  }-x_0-x_0\tau _\varepsilon 
\left(a+b\right)\right] -\frac{a}{b}. 
\end{align} 
Below we show that this is estimator of substitution. 
\begin{theorem}
\label{T3} Let  the condition  ${\scr A}$ be fulfilled and $\tau _\varepsilon
=\varepsilon ^\gamma $, where $ \gamma \in
\left(\frac{2}{5},\frac{2}{3}\right)$. Then the estimator $\bar\vartheta
_{\tau_\varepsilon }$ is consistent and asymptotically normal
\begin{align}
\label{17}
\frac{\tau _\varepsilon ^{3/2}}{\varepsilon }  \left(\bar\vartheta
_{\tau_\varepsilon }-\vartheta _0\right)\Longrightarrow  {\cal
  N}\left(0,D^2\right),\qquad D^2=\frac{4}{\left(a+b\right)^2b^2}. 
\end{align}
Moreover we have the convergence of moments too: for any $p>0$
\begin{align}
\label{18}
\varepsilon ^{-p\left(1-\frac{3}{2}\gamma \right)}\Ex_{\vartheta _0}\left|  \left(\bar\vartheta
_{\tau_\varepsilon }-\vartheta _0\right)\right|^p\longrightarrow
D^{\frac{p}{2}}\; \Ex \left|\xi \right|^p,\qquad \xi \sim {\cal N}\left(0,1\right).
\end{align}

\end{theorem}
{\bf Proof.}  The function $x_t\left(\vartheta
\right)$ can be expanded at the vicinity of the point $t=0 $
\begin{align*}
x_t\left(\vartheta \right)&=x_0+\int_{0}^{t}\left[ax_s+bx_{\vartheta
    s}\right]{\rm d}s\\
&=x_0+ x_0 t\left[a+b\right]+\frac12
x_0t^2\left[a+b\right]\left[a+b\vartheta \right] +O\left(t^3\right)
\end{align*}
and we can write 
\begin{align*}
\vartheta _\tau =\frac{2}{\left(a+b\right)b\tau ^2}\left[x_\tau
  \left(\vartheta \right)-x_0-x_0\tau \left(a+b\right)\right] -\frac{a}{b}+O\left(\tau \right).
\end{align*}
Remind that the process $X_t$ admits the representation
\begin{align*}
X_t=x_t\left(\vartheta _0\right)+\varepsilon
x_t^{\left(1\right)}\left(\vartheta _0\right){\rm d}t+O\left(\varepsilon ^2 \right)\quad 
\end{align*}
where the Gaussian process $x_t^{\left(1\right)}\left(\vartheta _0\right) $
(derivative of $X_t$ w.r.t. $\varepsilon $) 
has the stochastic differential
\begin{align*}
{\rm d}x_t^{\left(1\right)}\left(\vartheta
_0\right)=\left[ax_t^{\left(1\right)}\left(\vartheta
_0\right)+bx_{\vartheta _0t}^{\left(1\right)}\left(\vartheta _0\right)\right]
{\rm d}t +{\rm d}W_t,\quad x_0^{\left(1\right)}\left(\vartheta 
_0\right) =0
\end{align*}
and the expansion of it is
\begin{align*}
x_{\tau _\varepsilon }^{\left(1\right)}\left(\vartheta_0\right)=W_{\tau
  _\varepsilon }+O\left(\tau _\varepsilon ^2\right) .
\end{align*}
Hence
\begin{align*}
\bar\vartheta _{\tau_\varepsilon }& =\frac{2}{\left(a+b\right)b\tau
  _\varepsilon ^2}\left[x_{\tau _\varepsilon }\left(\vartheta
  _0\right)-x_0-x_0\tau _\varepsilon \left(a+b\right)\right] -\frac{a}{b}+
\frac{2 x_{\tau _\varepsilon }^{\left(1\right)}\left(\vartheta _0\right)
}{\left(a+b\right)b} \frac{\varepsilon}{\tau _\varepsilon ^2}\\
 &=\vartheta
_0+ \frac{2 W_{\tau _\varepsilon } }{\left(a+b\right)b}
\frac{\varepsilon}{\tau _\varepsilon ^2}+O\left(\tau _\varepsilon
\right)=\vartheta _0+ \frac{2 W_{1} }{\left(a+b\right)b}
\frac{\varepsilon}{\tau _\varepsilon ^{3/2}}+O\left(\tau _\varepsilon \right).
\end{align*}
Here $W_1 \sim {\cal N}\left(0,1\right)$.  Therefore 
\begin{align*}
\frac{\tau _\varepsilon ^{3/2}}{\varepsilon}\left(\bar\vartheta
_{\tau_\varepsilon }-\vartheta _0\right)= \frac{2 W_{1} }{\left(a+b\right)b}
+O\left(\frac{\tau _\varepsilon ^{5/2}}{\varepsilon}\right)
\end{align*}
and
if we suppose that 
\begin{align}
\label{19}
\frac{\tau _\varepsilon ^{3/2}}{\varepsilon}\longrightarrow \infty ,\qquad
\qquad \frac{\tau _\varepsilon ^ {5/2}}{\varepsilon }\longrightarrow 0 ,
\end{align}
then  we obtain the asymptotic normality \eqref{17}.
Note that the choice  $\tau _\varepsilon =\varepsilon ^\gamma $ with $\gamma \in
\left(\frac{2}{5},\frac{2}{3}\right)$ provides \eqref{19}.

The rate of convergence of this estimator is $\varepsilon ^\kappa $ with
$0<\kappa <\frac{2}{5}$, i.e.,
\begin{align}
\label{}
\frac{\bar\vartheta _{\tau_\varepsilon }-\vartheta _0}{\varepsilon ^\kappa
}\Longrightarrow {\cal N}\left(0,D^2\right). 
\end{align}
This rate of convergence is not sufficient for the construction of the
One-step MLE. We need an estimator with the rate $\kappa >\frac{1}{2}$. That
is why we propose below another estimator which uses  $\bar\vartheta
_{\tau_\varepsilon } $  and has   the rate 
of convergence up to the values grater than $\frac{1}{2}$.

\bigskip

{\bf One-step MDE.} Let us consider the equation \eqref{2.1}, where $X=X^{\tau
  _\varepsilon } =\left(X_t,0\leq t\leq \tau _\varepsilon \right)$. The
solution $\check\vartheta _{\tau _\varepsilon }$ of it satisfies the following
minimum distance equation
\begin{align}
\label{MDEq}
\int_{0}^{\tau _\varepsilon }\left[X_t-x_t\left(\check\vartheta _{\tau
    _\varepsilon }\right)\right] \dot x_t\left(\check\vartheta _{\tau
    _\varepsilon }\right){\rm d}t=0.
\end{align}
Here $\dot x_t\left(\vartheta \right)=\partial x_t\left(\vartheta
\right)/\partial \vartheta $ is solution of the equation
\begin{align}
\label{23}
\frac{{\rm d}\dot x_t\left(\vartheta \right)}{{\rm d}t}=a\dot
x_t\left(\vartheta \right)+ b\dot x_{\vartheta t}\left(\vartheta
\right)+bt\left[ax_{\vartheta t}\left(\vartheta \right)+bx_{\vartheta^2
    t}\left(\vartheta \right)\right] ,\qquad \dot x_{0}\left(\vartheta
\right)=0. 
\end{align}
 Let us write $x_t\left(\check\vartheta _{\tau _\varepsilon
 }\right)=x_t\left(\vartheta _{0 }\right)+\left(\check\vartheta _{\tau
   _\varepsilon }-\vartheta _0 \right)\dot x_t(\tilde\vartheta _{\tau
   _\varepsilon }) $, where the value $\tilde\vartheta _{\tau _\varepsilon }$
 satisfies the relation $|\tilde\vartheta _{\tau _\varepsilon }-\vartheta
 _0 |\leq |\check\vartheta _{\tau _\varepsilon }-\vartheta _0
 |$. From the equation \eqref{MDEq} we obtain the representation
\begin{align*}
\check\vartheta _{\tau _\varepsilon }=\vartheta _0+\left(\int_{0}^{\tau _\varepsilon
}\dot x_t(\tilde\vartheta _{\tau _\varepsilon })\dot x_t(\check\vartheta
_{\tau _\varepsilon }){\rm d}t\right)^{-1} \int_{0}^{\tau _\varepsilon
}\left[X_t-x_t\left(\vartheta _{0}\right)\right] \dot x_t\left(\check\vartheta _{\tau
    _\varepsilon }\right){\rm d}t.
\end{align*}
Let us denote
\begin{align*}
Q_{\tau _\varepsilon }\left(\vartheta \right)=\int_{0}^{\tau _\varepsilon}\dot x_t(\vartheta)^2{\rm d}t. 
\end{align*}
The One-step MDE $\vartheta _{\tau _\varepsilon }^*$  we introduce as follows
\begin{align}
\label{}
\vartheta _{\tau _\varepsilon }^*=\bar\vartheta _{\tau_\varepsilon }+Q_{\tau
  _\varepsilon }\left(\bar\vartheta _{\tau_\varepsilon } \right)^{-1} \int_{0}^{\tau _\varepsilon
}\left[X_t-x_t\left(\bar\vartheta _{\tau_\varepsilon }\right)\right] \dot
x_t\left(\bar\vartheta _{\tau_\varepsilon }\right){\rm d}t.
\end{align}

\begin{theorem}
\label{T4}  Let the conditions of Theorem \ref{T3} be fulfilled, then for any
$p>0$ there exists a constant $C>0$ such that
\begin{align*}
\varepsilon ^{-p\delta }\Ex_{\vartheta _0}\left|\vartheta _{\tau _\varepsilon
}^*-\vartheta _0\right|^p \leq C
\end{align*}
where $\delta \in \left(0,\frac{2}{3}\right)$.
\end{theorem}
{\bf Proof.} Let us write
\begin{align*}
\frac{\vartheta _{\tau _\varepsilon }^*-\vartheta _0}{\varepsilon ^\delta
}=\frac{\bar\vartheta _{\tau_\varepsilon }-\vartheta _0}{\varepsilon ^\delta } +\varepsilon ^{-\delta }Q_{\tau
  _\varepsilon }\left(\bar\vartheta _{\tau_\varepsilon } \right)^{-1} \int_{0}^{\tau _\varepsilon
}\left[X_t-x_t\left(\bar\vartheta _{\tau_\varepsilon }\right)\right] \dot
x_t\left(\bar\vartheta _{\tau_\varepsilon }\right){\rm d}t.
\end{align*}
We have
\begin{align*}
\varepsilon ^{-\delta }\left[X_t-x_t\left(\bar\vartheta _{\tau_\varepsilon
  }\right)\right]=\varepsilon ^{-\delta }\left[X_t-x_t\left(\vartheta
  _0\right)\right]-\varepsilon ^{-\delta } \left(\bar\vartheta
_{\tau_\varepsilon } -\vartheta _0\right)\dot x_t(\tilde\vartheta _{\tau
  _\varepsilon }).
\end{align*}
Recall that
\begin{align*}
\Ex_{\vartheta _0}\left|X_t-x_t\left(\vartheta _0\right)\right|^p\leq
C\,\varepsilon ^p 
\end{align*}
(see \eqref{A5} below). 
Therefore for any $\delta <1$ we have $\varepsilon ^{-\delta
}\left[X_t-x_t\left(\vartheta _0\right)\right]=o\left(1\right) $. Further, we
can write once more the expansion $\dot x_t(\tilde\vartheta _{\tau
  _\varepsilon })=\dot x_t(\bar\vartheta _{\tau
  _\varepsilon })+(\tilde\vartheta _{\tau
  _\varepsilon }-\bar\vartheta _{\tau
  _\varepsilon })\ddot x_t(\tilde{\bar\vartheta }_{\tau
  _\varepsilon }) $. Hence 
\begin{align*}
&\varepsilon ^{-\delta }Q_{\tau
  _\varepsilon }\left(\bar\vartheta _{\tau_\varepsilon } \right)^{-1} \int_{0}^{\tau _\varepsilon
}\left[X_t-x_t\left(\bar\vartheta _{\tau_\varepsilon }\right)\right] \dot
x_t\left(\bar\vartheta _{\tau_\varepsilon }\right){\rm d}t\\
& \qquad =o\left(1\right)-\varepsilon ^{-\delta } \left(\bar\vartheta
_{\tau_\varepsilon } -\vartheta _0\right)Q_{\tau
  _\varepsilon }\left(\bar\vartheta _{\tau_\varepsilon } \right)^{-1} \int_{0}^{\tau _\varepsilon
}\dot x_t\left(\bar\vartheta _{\tau_\varepsilon }\right)^2{\rm d}t\\
& \qquad \qquad +\varepsilon ^{-\delta } \left(\bar\vartheta
_{\tau_\varepsilon } -\vartheta _0\right) (\bar\vartheta
_{\tau_\varepsilon } -\tilde\vartheta_{\tau_\varepsilon })Q_{\tau
  _\varepsilon }\left(\bar\vartheta _{\tau_\varepsilon } \right)^{-1} \int_{0}^{\tau _\varepsilon
}\ddot x_t(\tilde{\bar\vartheta} _{\tau_\varepsilon })\dot
x_t\left(\bar\vartheta _{\tau_\varepsilon }\right){\rm d}t\\
& \qquad =o\left(1\right)- \frac{\left(\bar\vartheta
_{\tau_\varepsilon } -\vartheta _0\right)}{\varepsilon ^\delta }+ \frac{\left(\bar\vartheta
_{\tau_\varepsilon } -\vartheta _0\right) (\bar\vartheta
_{\tau_\varepsilon } -\tilde\vartheta_{\tau_\varepsilon })}{\varepsilon
  ^\delta }\frac{ \int_{0}^{\tau _\varepsilon 
}\ddot x_t(\tilde{\bar\vartheta} _{\tau_\varepsilon })\dot
x_t\left(\bar\vartheta _{\tau_\varepsilon }\right) {\rm d}t}{ Q_{\tau
  _\varepsilon }\left(\bar\vartheta _{\tau_\varepsilon } \right)}
\end{align*}
For the small values of $t$ the solution of the equation \eqref{23}  can be
written as follows
\begin{align*}
\dot x_t\left(\vartheta \right)=\frac{1}{2}bt^2x_0\left(a+B\right)+O\left(t^3\right).
\end{align*}
For the second derivative we have the similar expansions
\begin{align*}
\ddot x_t\left(\vartheta
\right)=\frac{1}{2}at^2x_0\left(a+b\right)+\frac{1}{2}bt^2x_0\left(a+b\right)
+O\left(t^3\right) 
\end{align*}
Hence we have estimates
\begin{align*}
Q_{\tau
  _\varepsilon }\left(\bar\vartheta _{\tau_\varepsilon } \right)\geq c\tau
_\varepsilon ^5,\qquad  \quad \int_{0}^{\tau _\varepsilon 
}\ddot x_t(\tilde{\bar\vartheta} _{\tau_\varepsilon })\dot
x_t\left(\bar\vartheta _{\tau_\varepsilon }\right) {\rm d}t\leq C\tau
_\varepsilon ^5.
\end{align*}
Note as well that
\begin{align*}
\left|\left(\bar\vartheta
_{\tau_\varepsilon } -\vartheta _0\right) (\bar\vartheta
_{\tau_\varepsilon } -\tilde\vartheta_{\tau_\varepsilon })\right|\leq 2\left|\bar\vartheta
_{\tau_\varepsilon } -\vartheta _0\right|^2.
\end{align*}
Therefore we obtain the relation
\begin{align*}
\left|\frac{\left(\bar\vartheta
_{\tau_\varepsilon } -\vartheta _0\right) (\bar\vartheta
_{\tau_\varepsilon } -\tilde\vartheta_{\tau_\varepsilon })}{\varepsilon
  ^\delta }\frac{ \int_{0}^{\tau _\varepsilon 
}\ddot x_t(\tilde{\bar\vartheta} _{\tau_\varepsilon })\dot
x_t\left(\bar\vartheta _{\tau_\varepsilon }\right) {\rm d}t}{ Q_{\tau
  _\varepsilon }\left(\bar\vartheta _{\tau_\varepsilon } \right)}\right|\leq C
\frac{\left|\bar\vartheta
_{\tau_\varepsilon } -\vartheta _0\right|^2 }{\varepsilon ^\delta
}=\frac{\varepsilon ^{2\kappa }}{\varepsilon ^\delta } \left|\bar
u_\varepsilon \right|^2, 
\end{align*}
where $\bar u_\varepsilon=\varepsilon ^{-\kappa }\left( \bar\vartheta
_{\tau_\varepsilon } -\vartheta _0\right)$. If we take $\delta <2\kappa $,
i.e., $0<\delta <\frac{2}{3}$, then for the One-step MDE we obtain
\begin{align*}
\frac{\vartheta _{\tau _\varepsilon }^*-\vartheta
  _0}{\varepsilon ^\delta } =o\left(1\right),\qquad \quad   \Ex_{\vartheta
  _0}\left|\frac{\vartheta _{\tau _\varepsilon }^*-\vartheta 
  _0}{\varepsilon ^\delta }\right|^p=o\left(1\right).
\end{align*}
Now the rate of convergence of preliminary estimator $\vartheta _{\tau _\varepsilon }^* $ can be grater than $
\frac12$ and it can be used for the construction of One-step MLE.

\bigskip

{\bf Two-step MLE.}
Having this preliminary estimator we construct the Two-step MLE as follows
\begin{align}
\label{21}
\vartheta _\varepsilon ^\star=\vartheta^* _{\tau _\varepsilon }+{\rm
  I}_{\tau _\varepsilon }\left(\vartheta^* _{\tau _\varepsilon
}\right)^{-1}b \int_{\tau _\varepsilon }^{T}t\left[ax_{t}+bx_{\vartheta^*
    _{\tau _\varepsilon }t}\right]\left[{\rm
    d}X_t-\left(aX_t+bX_{\vartheta^* _{\tau _\varepsilon }t}\right){\rm
    d}t\right],
\end{align}
where
\begin{align*}
{\rm I}_\tau \left(\vartheta \right)=b^2 \int_{\tau }^{T}t^2 \left[ax_{\vartheta   
  t}\left(\vartheta \right)+bx_{\vartheta^2
  t}\left(\vartheta \right)\right]^2 {\rm d}t.
\end{align*}
Note that we have no more problem with definition of the stochastic integral
because the estimator $\vartheta^* _{\tau _\varepsilon }$ depends on the
observations $X^{\tau _\varepsilon }=\left(X_t,0\leq t\leq \tau _\varepsilon \right)$
 and the integral starts at the moment $\tau _\varepsilon $.

\begin{theorem}
\label{T5} Suppose that the condition ${\scr A}$ is fulfilled  and $\tau _\varepsilon
=\varepsilon ^\gamma $ with $\gamma \in
\left(\frac{2}{5},\frac{1}{2}\right)$. Then the estimator $\vartheta
_\varepsilon ^\star$ is consistent, asymptotically normal 
\begin{align}
\label{}
\varepsilon ^{-1}\left(\vartheta _\varepsilon ^\star-\vartheta
_0\right)\Longrightarrow {\cal N}\left(0, {\rm
  I}\left(\vartheta _{0}\right)^{-1}\right)
\end{align}
and we have the convergence of polynomial moments. 
\end{theorem}
{\bf Proof.} Below we substitute the observations \eqref{14}, where $\vartheta
=\vartheta _0$, in the stochastic integral 
\begin{align*}
\varepsilon ^{-1}\left(\vartheta _\varepsilon ^\star-\vartheta
_0\right)&=\varepsilon ^{-1}\left(\vartheta^* _{\tau _\varepsilon
}-\vartheta _0\right)+{\rm 
  I}_{\tau _\varepsilon }\left(\vartheta^* _{\tau _\varepsilon }\right)^{-1} b\int_{\tau
  _\varepsilon }^{T}t\left[ax_{t}+bx_{\vartheta^* _{\tau _\varepsilon
    }t}\right]{\rm d}W_t\\
&\qquad +\varepsilon ^{-1} {\rm
  I}_{\tau _\varepsilon }\left(\vartheta^* _{\tau _\varepsilon }\right)^{-1}b^2 \int_{\tau
  _\varepsilon }^{T}t\left[ax_{t}+bx_{\vartheta^* _{\tau _\varepsilon
    }t}\right]\left(X_{\vartheta _0t}-X_{\vartheta^* _{\tau
      _\varepsilon }t}\right){\rm d}t.
\end{align*}
From the consistency of the estimator $\vartheta^* _{\tau _\varepsilon }$
and continuity of the corresponding functions we
obtain the convergences
\begin{align*}
 b^2\int_{\tau
  _\varepsilon }^{T}t^2\left[ax_{t}+bx_{\vartheta^* _{\tau _\varepsilon
    }t}\right]^2{\rm d}t\longrightarrow {\rm
  I}\left(\vartheta _{0}\right),\qquad {\rm
  I}_{\tau _\varepsilon }\left(\vartheta^* _{\tau _\varepsilon }\right)\longrightarrow {\rm
  I}\left(\vartheta _{0}\right).
\end{align*}
Hence by the central limit theorem for stochastic integrals we have the
asymptotic normality 
\begin{align*}
{\rm
  I}_{\tau _\varepsilon }\left(\vartheta^* _{\tau _\varepsilon }\right)^{-1} b\int_{\tau
  _\varepsilon }^{T}t\left[ax_{t}+bx_{\vartheta^* _{\tau _\varepsilon
    }t}\right]{\rm d}W_t\Longrightarrow {\cal N}\left(0, {\rm
  I}\left(\vartheta _{0}\right)^{-1}\right).
\end{align*}
Remark, that for this model of observations we have the asymptotic normality
``in probability'' too
\begin{align*}
{\rm
  I}_{\tau _\varepsilon }\left(\vartheta^* _{\tau _\varepsilon }\right)^{-1}b \int_{\tau
  _\varepsilon }^{T}t\left[ax_{t}+bx_{\vartheta^* _{\tau _\varepsilon
    }t}\right]{\rm d}W_t\longrightarrow {\rm
  I}\left(\vartheta _{0}\right)^{-1} b \int_{0 }^{T}t\left[ax_{t}+bx_{\vartheta _{0
    }t}\right]{\rm d}W_t,
\end{align*}
where
\begin{align*}
{\rm
  I}\left(\vartheta _{0}\right)^{-1} b \int_{0 }^{T}t\left[ax_{t}+bx_{\vartheta _{0
    }t}\right]{\rm d}W_t\sim {\cal N}\left(0, {\rm
  I}\left(\vartheta _{0}\right)^{-1}\right).
\end{align*}
Further, uniformly on $t\in \left[0,T\right]$
\begin{align*}
&X_{\vartheta^* _{\tau _\varepsilon }t}-X_{\vartheta _0t}=x_{\vartheta^*
  _{\tau _\varepsilon }t}\left(\vartheta _0\right)-x_{\vartheta
  _0t}\left(\vartheta _0\right)+\varepsilon \left[x_{\vartheta^* _{\tau
      _\varepsilon }t}^{\left(1\right)}\left(\vartheta _0\right)-x_{\vartheta
    _0t}^{\left(1\right)}\left(\vartheta _0\right) \right]+O\left(\varepsilon
^2\right)\\ 
&\qquad \qquad =\left(\vartheta^* _{\tau _\varepsilon } -\vartheta
_0\right)t\left[ax_{\vartheta _0 t}+bx_{\vartheta _0
    ^2t}\right]+\varepsilon\Delta _\varepsilon  +O\left(\varepsilon 
^2\right)+O\left(\left(\vartheta^* _{\tau _\varepsilon } -\vartheta
_0 \right)^2\right),
\end{align*}
where
\begin{align*}
\Delta _\varepsilon &=x_{\vartheta^* _{\tau
      _\varepsilon }t}^{\left(1\right)}\left(\vartheta _0\right)-x_{\vartheta
    _0t}^{\left(1\right)}\left(\vartheta _0\right)\\
& =\int_{\vartheta _0t}^{\vartheta^* _{\tau
      _\varepsilon }t}\left[ ax_s^{\left(1\right)}\left(\vartheta
  _0\right)+bx_{\vartheta _0s}^{\left(1\right)}\left(\vartheta
  _0\right)\right]{\rm d}s +W_{\vartheta^* _{\tau
      _\varepsilon }t}-W_{\vartheta _0t}\\
& =O\left(\vartheta^* _{\tau _\varepsilon } -\vartheta
_0\right)+O\left(\left|\vartheta^* _{\tau _\varepsilon } -\vartheta
_0\right|^{1/2}\right)=O\left(\left|\vartheta^* _{\tau _\varepsilon } -\vartheta
_0\right|^{1/2}\right).
\end{align*}
Therefore
\begin{align*}
X_{\vartheta^* _{\tau _\varepsilon }t}-X_{\vartheta _0t}=\left(\vartheta^*
_{\tau _\varepsilon } -\vartheta _0\right)t\left[ax_{\vartheta _0
    t}+bx_{\vartheta _0 ^2t}\right]+\varepsilon \left|\vartheta^* _{\tau _\varepsilon } -\vartheta
_0\right|^{1/2}O\left(1\right)
\end{align*}
and we can write
\begin{align*}
&\varepsilon ^{-1} {\rm I}_{\tau _\varepsilon }\left(\vartheta^* _{\tau _\varepsilon
  }\right)^{-1}b^2 \int_{\tau _\varepsilon
  }^{T}t\left[ax_{t}+bx_{\vartheta^* _{\tau _\varepsilon
      }t}\right]\left(X_{\vartheta _0t}-X_{\vartheta^* _{\tau _\varepsilon
    }t}\right){\rm d}t\\ 
&\quad =-\frac{\left(\vartheta^* _{\tau
    _\varepsilon } -\vartheta _0\right)}{\varepsilon } {\rm I}_{\tau _\varepsilon }\left(\vartheta^* _{\tau
    _\varepsilon }\right)^{-1}b^2\int_{\tau _\varepsilon
  }^{T}t^2\left[ax_{t}+bx_{\vartheta^* _{\tau _\varepsilon }t}\right]
  \left[ax_{t}+bx_{\vartheta _0t}\right] {\rm d}t\\
&\qquad \qquad \qquad +\left|\vartheta^* _{\tau
    _\varepsilon } -\vartheta _0\right|^{1/2}O\left(1\right)\\
&\quad =-\frac{\left(\vartheta^* _{\tau
    _\varepsilon } -\vartheta _0\right)}{\varepsilon }\left[1 +O\left( \vartheta^* _{\tau
    _\varepsilon } -\vartheta _0\right)
  \right]+\left|\vartheta^* _{\tau
    _\varepsilon } -\vartheta _0\right|^{1/2}O\left(1\right).
\end{align*}
Here by conditions of the Theorem \ref{T5} 
\begin{align*}
\frac{\left(\vartheta^* _{\tau _\varepsilon } -\vartheta _0\right)^2
}{\varepsilon }=\frac{\varepsilon ^{2\delta }}{\varepsilon  }O\left(1\right)=o\left(1\right)
\end{align*}
Finally we have
\begin{align*}
\varepsilon ^{-1}{\left(\vartheta _\varepsilon ^\star-\vartheta _0\right)}{}&=
\frac{b }{{\rm 
  I}\left(\vartheta _{0}\right) }\int_{0 }^{T}t\left[ax_{t}+bx_{\vartheta _{0
    }t}\right]{\rm d}W_t+o\left(1\right). 
\end{align*}
It can be verified that we have the uniform on compacts $\KK\subset\Theta $
convergence of moments too
\begin{align*}
\varepsilon ^{-2}\sup_{\vartheta _0\in \KK}\Ex_{\vartheta _0}\left|\vartheta
_\varepsilon ^\star-\vartheta _0\right|^2\longrightarrow \sup_{\vartheta _0\in
  \KK}{\rm I}\left(\vartheta _{0}\right).
\end{align*}
Therefore for any (small) $\nu >0$
\begin{align*}
\lim_{\varepsilon \rightarrow 0}\varepsilon ^{-2}\sup_{\left|\vartheta
  -\vartheta _0\right|\leq \nu }\Ex_{\vartheta }\left|\vartheta _\varepsilon
^\star-\vartheta\right|^2=\sup_{\left|\vartheta
  -\vartheta _0\right|\leq \nu}{\rm I}\left(\vartheta \right)
\end{align*}
and
\begin{align*}
\lim_{\nu \rightarrow 0}\sup_{\left|\vartheta
  -\vartheta _0\right|\leq \nu}{\rm I}\left(\vartheta \right)={\rm
  I}\left(\vartheta _{0}\right). 
\end{align*}
Hence Two-step MLE  $ \vartheta _\varepsilon
^\star$ is asymptotically efficient (see \eqref{HLC}).

\bigskip

The procedure of estimation is the following: first we calculate the estimator
of substitution $\bar\vartheta _{\tau_\varepsilon }$, then with the help of
this estimator we construct the Two-step MDE $\vartheta _{\tau _\varepsilon
}^* $, which has the ``good rate'' and can be used for construction of the
Two-step MLE $\vartheta _\varepsilon ^\star $. This last estimator is
asymptotically equivalent to the MLE.

Let us write together all estimators. The first one is estimator of
substitution 
\begin{align*}
\bar\vartheta _{\tau_\varepsilon }
=\frac{2}{\left(a+b\right)b\tau_{\varepsilon } ^2}\left[X_{\tau _\varepsilon
  }-x_0-x_0\tau _\varepsilon 
\left(a+b\right)\right] -\frac{a}{b}.
\end{align*}
Then, using this estimator we calculate the Two-step MDE
\begin{align*}
\vartheta _{\tau _\varepsilon }^*=\bar\vartheta _{\tau_\varepsilon }+Q_{\tau
  _\varepsilon }\left(\bar\vartheta _{\tau_\varepsilon } \right)^{-1} \int_{0}^{\tau _\varepsilon
}\left[X_t-x_t\left(\bar\vartheta _{\tau_\varepsilon }\right)\right] \dot
x_t\left(\bar\vartheta _{\tau_\varepsilon }\right){\rm d}t.
\end{align*}
Finally, we can calculate the Two-step MLE
\begin{align*}
\vartheta _\varepsilon ^\star=\vartheta^* _{\tau _\varepsilon }+{\rm
  I}_{\tau _\varepsilon }\left(\vartheta^* _{\tau _\varepsilon
}\right)^{-1}b \int_{\tau _\varepsilon }^{T}t\left[ax_{t}+bx_{\vartheta^*
    _{\tau _\varepsilon }t}\right]\left[{\rm
    d}X_t-\left(aX_t+bX_{\vartheta^* _{\tau _\varepsilon }t}\right){\rm
    d}t\right].
\end{align*}
The advantage of this procedure is the absence of the problems of optimization
like \eqref{mle} and  \eqref{2.1}.

\bigskip

{\bf Remark.} Note that it is possible to introduce the {\it Two-step MLE-process}
$\vartheta _{t,\varepsilon} ^\star , \tau _\varepsilon <t\leq T$, which can be
used in the problems, where we need estimator for all $t\in (0,T]$. For
  example, such estimators are used in the approximation of the solutions of
  backward stochastic differential equations \cite{Kut14}. Here this
  estimator-process can be written as follows
\begin{align*}
\vartheta _{t,\varepsilon} ^\star=\vartheta^* _{\tau _\varepsilon }+{\rm
  I}^t_{\tau _\varepsilon }\left(\vartheta^* _{\tau _\varepsilon
}\right)^{-1}b \int_{\tau _\varepsilon }^{t}s\left[ax_{s}+bx_{\vartheta^*
    _{\tau _\varepsilon }s}\right]\left[{\rm
    d}X_s-\left(aX_s+bX_{\vartheta^* _{\tau _\varepsilon }s}\right){\rm
    d}s\right],
\end{align*}
where $t\in(\tau _\varepsilon, T] $ and
\begin{align*}
{\rm I}^t_{\tau }\left(\vartheta\right)=b^2\int_{\tau }^{t}\dot
x_s\left(\vartheta \right)^2{\rm d}s .
\end{align*}
It can be shown that 
\begin{align*}
\varepsilon ^{-1}\left(\vartheta _{t,\varepsilon} ^\star-\vartheta
_0\right)\Longrightarrow {\cal N}\left(0, {\rm I}^t_{\tau
}\left(\vartheta_0\right)^{-1}\right).
\end{align*}
For the proof of this convergence we can use the proof of the Theorem
\ref{T5}, where we just put $T=t$. Morover, it can be shown  that we
have the weak convergence of the stochastic process $\eta _{t,\varepsilon
}={\rm I}^t_{\tau }\left(\vartheta\right)$ (see \cite{KZ14}).

\subsection{Nonlinear scale-type  equation}

We consider the following  nonlinear SDE 
\begin{align}
\label{28}
{\rm d}X_t=S\left(X_{\vartheta t}\right)\,{\rm d}t+\varepsilon {\rm d}W_t,
\quad X_0=x_0,\quad 0\leq t\leq T, 
\end{align}
where $\vartheta \in\Theta =\left(\alpha ,\beta \right)$, $0<\alpha <\beta
<1$. Suppose that the function $S\left(x\right)$ has four continuous bounded
derivatives. Recall that this equation has a unique strong solution (see
Theorem 4.6 in \cite{LS01}).

The limit differential equation is
\begin{align*}
\frac{{\rm d}x_t}{{\rm d}t}=S\left(x_{\vartheta t}\right),\qquad x_0,\quad 0\leq t\leq T.
\end{align*}
Its solution $x_t=x_t\left(\vartheta \right)$ is a function of $\vartheta $. 

The preliminary MDE $\bar\vartheta _\varepsilon $ for this model of
observations is defined by the same equation as in \eqref{2.1}. This estimator
is consistent, asymptotically normal and the moments converge. The proof
is similar to  the proof of Theorem \ref{T1}.

The Fisher information is
\begin{align*}
{\rm I}\left(\vartheta \right)=\int_{0}^{T}t^2S'\left(x_{\vartheta
  t}\left(\vartheta \right)\right)^2 S\left(x_{\vartheta^2
  t}\left(\vartheta \right)\right)^2 {\rm d}t.
\end{align*}
The properties of the MLE were studied in \cite{A86}, \cite{Kut88} (linear
case) and \cite{Kut94}. It was shown that the MLE $\hat\vartheta _\varepsilon
$ is consistent, asymptotically normal
\begin{align*}
\varepsilon ^{-1}\left(\hat\vartheta _\varepsilon -\vartheta
_0\right)\Longrightarrow  {\cal N}\left(0,{\rm I}
\left(\vartheta _0\right)^{-1}\right),
\end{align*}
 and asymptotically efficient. Its numerical calculation can have the same
 difficulties as in the case of observations \eqref{1}.

 To introduce the One-step MLE we
recall the definition of the functions $H\left(\cdot \right)$ and $\Psi
\left(\cdot \right)$ for this model:
\begin{align*}
H\left(t,\vartheta \right)&=tS'\left(x_{\vartheta
  t}\left(\vartheta \right)\right) S\left(x_{\vartheta^2
  t}\left(\vartheta \right)\right),\\
\Psi \left(\vartheta \right)&=H\left(T,\vartheta \right)X_T-H\left(0,\vartheta
\right)x_0- \int_{0}^{T}H'_t\left(t,\vartheta \right)X_t\,{\rm d}t.
\end{align*}
Note that the function $H\left(\cdot \right)$ and its derivatives are bounded
on the set $\left[0,T\right]\times \Theta $.
Then we can define
\begin{align*}
\vartheta _\varepsilon ^\star&=\bar\vartheta _\varepsilon +{\rm I}
\left(\bar\vartheta _\varepsilon\right)^{-1} \left[\Psi
  \left(\bar\vartheta _\varepsilon \right)-
  \int_{0}^{T}t\,S'\left(x_{\bar\vartheta _\varepsilon t}\left(\bar\vartheta
  _\varepsilon \right)\right) S\left(x_{t\bar\vartheta
    _\varepsilon^2}\left(\bar\vartheta _\varepsilon \right)\right)
  S\left(X_{\bar\vartheta _\varepsilon t}\right){\rm d}t \right] .
\end{align*}

This One-step MLE is consistent and asymptotically normal
\begin{align*}
\varepsilon ^{-1}\left(\vartheta _\varepsilon ^\star-\vartheta
_0\right)\Longrightarrow  {\cal N}\left(0,{\rm I}
\left(\vartheta _0\right)^{-1}\right).
\end{align*}
The proof follows the same lines as the proof of the Theorem \ref{T2}.

\section{Discussion}

The proposed result can be generalized on the slightly more  general
models. For example, the cases  with observations
\begin{align*}
{\rm d}X_t=\left[S_1\left(X_{t} \right)+S_2\left(X_{t-\vartheta } \right)\right]{\rm
  d}t +\varepsilon 
\sigma \left(t,X_t \right){\rm d}W_t,\;  X_s=x_s, s\leq 0,\; 0\leq t\leq T
\end{align*}
and
\begin{align*}
{\rm d}X_t=S\left(\vartheta ,t,X_{t-\vartheta} \right){\rm d}t+\varepsilon
\sigma \left(t,X_t \right){\rm d}W_t,\quad  X_s=x_s, s\leq 0,\quad 0\leq t\leq T
\end{align*}
can be treated exactly as it was done above. Here the functions
$S_1\left(\cdot \right), S_2\left(\cdot \right),$ $ S\left(\cdot
\right),\sigma \left(\cdot \right)$ are supposed to be deterministic and
smooth.

Another model, which can be treated by a similar way is
\begin{align*}
{\rm d}X_t=S\left(X_{t-f\left(\vartheta ,t\right)} \right){\rm
  d}t +\varepsilon 
{\rm d}W_t,\;  X_s=x_s, s\leq 0,\; 0\leq t\leq T.
\end{align*}
Here $f\left(\vartheta ,t\right)\geq 0$ is some known smooth function. If
$f\left(\vartheta ,t\right)=\vartheta $ and $f\left(\vartheta
,t\right)=t-\vartheta t  $ then we obtain the models \eqref{1} and \eqref{28}
respectively.

It is possible to study the models with multiple delays like
\begin{align*}
{\rm d}X_t=\sum_{k=1}^{K}S_k\left(X_{t -\vartheta_k} \right){\rm d}t+\varepsilon
{\rm d}W_t,\quad  X_s=x_s, s\leq 0,\quad 0\leq t\leq T.
\end{align*} 
In the last model the delays form a vector $\vartheta =\left(\vartheta
_1,\ldots,\vartheta _K\right)$, where $\vartheta
_1<\vartheta _2<\ldots<\vartheta _K\ $ and the problem of estimation $\vartheta $
became more complicate but the construction of the One-step procedure can be
realized using the corresponding modification of the given above procedure.

\section{Appendix}\label{app}

We suppose that the function $S\left(x\right),  x\in
{\cal R}$ has two continuous bounded
derivatives  $\left|S'\left(x\right)\right|\leq L$,
$\left|S''\left(x\right)\right|\leq M$, therefore  the equation \eqref{1} has a unique strong
solution (see Theorem 4.6 in \cite{LS01}). Moreover, we have the following
estimates (with probability 1)
\begin{align*}
&\left|X_t-x_t\right|\leq \int_{0}^{t}\left|S\left(X_{s-\vartheta _0 
}\right)-S\left(x_{s-\vartheta _0 }\right) \right|{\rm d}s+\varepsilon \left|W_t\right| \\
&\qquad \leq L\int_{0}^{t}\left|X_{s-\vartheta _0 
}-x_{s-\vartheta _0  } \right|{\rm d}s+\varepsilon \left|W_t\right|\leq L\int_{0}^{t}\left|X_{s
}-x_{s } \right|{\rm d}s+\varepsilon \sup_{0\leq s\leq T}\left|W_s\right|.
\end{align*}

Recall the well-known  lemma.
\begin{lemma}
\label{L2} {\rm [Gronwall-Bellman]}
Suppose that the function $V\left(t\right)\geq 0, 0\leq t\leq T$ satisfies the
inequality 
\begin{align*}
V\left(t\right)\leq A\int_{0}^{t}V\left(s\right){\rm d}s+B
\end{align*}
where the constants $A>0, B>0$. Then
\begin{align*}
V\left(t\right)\leq B\,e^{AT}.
\end{align*}
\end{lemma}

By this lemma we can write
\begin{align}
\label{A1}
\left|X_t-x_t\right|\leq \varepsilon\, C\, \widehat{W}_T,
\end{align}
where we denoted $\widehat{W}_T=\sup_{0\leq s\leq T}\left|W_s\right|$. 

For the second moment we have
\begin{align*}
\Ex_S\left|X_t-x_t\right|^2&\leq 2L^2T\int_{0}^{t}\Ex_S\left|X_{s-\vartheta _0 
}-x_{s-\vartheta _0  }\right|^2{\rm d}s+2\varepsilon ^2T\\
&\leq 2L^2T\int_{0}^{t}\Ex_S\left|X_{s
}-x_{s }\right|^2{\rm d}s+2\varepsilon ^2T
\end{align*}
and by the same lemma
\begin{align}
\label{A2}
\Ex_S\left|X_t-x_t\right|^2\leq C\;\varepsilon ^2.
\end{align}
Recall that for $s\leq 0$ we have $X_s-x_s=x_0-x_0=0$. 

 Let us write the formal derivative $X_t^{\left(1\right)}$  of $X_t$
 w.r.t. $\varepsilon $, then we obtain the equation
\begin{align*}
{\rm d}X_t^{\left(1\right)}=S'_x\left(X_{t-\vartheta _0
}\right)X_{t-\vartheta _0 }^{\left(1\right)}{\rm d}t+{\rm d}W_t,\qquad
X_s^{\left(1\right)}=0,\;{\rm for}\; s\leq 0,\quad 0\leq t\leq T.
\end{align*}
The Gaussian process $X_t^{\left(1\right)}|_{\varepsilon
  =0}=x_t^{\left(1\right)}\left(\vartheta _0\right)$ satisfies the linear equation
\begin{align}
\label{A4}
{\rm d}x_t^{\left(1\right)}=S'_x\left(x_{t-\vartheta _0
}\right)x_{t-\vartheta _0 }^{\left(1\right)}{\rm d}t+{\rm d}W_t,\quad
x_s^{\left(1\right)}=0,\;{\rm for}\; s\leq 0,\quad 0\leq t\leq T.
\end{align}
The proof of 
\begin{align*}
V_\varepsilon\left(t\right) =\Ex_{\vartheta
  _0}\left|\frac{X_t-x_t\left(\vartheta _0\right)}{\varepsilon
}-x_t^{\left(1\right)}\left(\vartheta \right) \right|^2\rightarrow 0
\end{align*}
can be done  using the standard technique based on
Gronwall-Bellman lemma as follows. We write
\begin{align*}
&v_\varepsilon \left(t\right)=\frac{X_t-x_t\left(\vartheta
    _0\right)}{\varepsilon
  }-x_t^{\left(1\right)}\\
 &\qquad=\int_{0}^{t}\left[\frac{S\left(X_{s-\vartheta
        _0}\right)-S\left(x_{s-\vartheta _0} \right)}{\varepsilon
    }-S'\left(x_{s-\vartheta _0} \right)x_t^{\left(1\right)} \right]{\rm
    d}s\\ 
&\qquad=\int_{0}^{t}\left[S'(\tilde X_{s-\vartheta
      _0})\frac{X_{s-\vartheta _0}-x_{s-\vartheta _0} }{\varepsilon
    }-S'\left(x_{s-\vartheta _0} \right)x_t^{\left(1\right)} \right]{\rm
    d}s\\
 &\qquad=\int_{0}^{t}S'(\tilde X_{s-\vartheta
    _0})\left[\frac{X_{s-\vartheta _0}-x_{s-\vartheta _0} }{\varepsilon
    }-S'\left(x_{s-\vartheta _0} \right)x_t^{\left(1\right)} \right]{\rm
    d}s\\
 &\qquad \qquad +\int_{0}^{t}\left[S'(\tilde X_{s-\vartheta
      _0})-S'\left(x_{s-\vartheta _0} \right)\right]{\rm d}s,
\end{align*}
where $\left|\tilde X_{s-\vartheta
    _0}-x_{s-\vartheta _0}\right|\leq \left| X_{s-\vartheta
    _0}-x_{s-\vartheta _0}\right|$.   Hence
\begin{align*}
\Ex_{\vartheta _0}\left|v_\varepsilon \left(t\right)\right|^2&\leq
L^2T\int_{0}^{t}\Ex_{\vartheta _0}\left|v_\varepsilon 
\left(s-\vartheta _0\right)\right|^2 {\rm d}s\\
&\qquad \qquad +M^2T\int_{0}^{T}\Ex_{\vartheta
  _0}\left| X_{t-\vartheta _0}-x_{t-\vartheta _0}\right|^2{\rm d}s 
\end{align*}
and
\begin{align*}
V_\varepsilon \left(t\right)&\leq L^2T\int_{0}^{t}V_\varepsilon \left(s-\vartheta _0\right)
{\rm d}s +M^2T^2C\varepsilon ^2 \\
& \leq L^2T\int_{0}^{t}V_\varepsilon \left(s\right)
{\rm d}s +M^2T^2C\varepsilon ^2  .
\end{align*}
 Therefore by Gronwall-Bellman lemma $V_\varepsilon \left(t\right)\rightarrow
0$.

\bigskip

Let us study the quantity $R_\varepsilon \left(t\right)=\Ex_{\vartheta
  _0}\left|X_t-x_t\left(\vartheta _0\right)\right|^p$, where
$x_t\left(\vartheta _0\right)$ and $X_t$ are solutions of the equations
\eqref{14a} and \eqref{14} respectively.  For the difference
$r_\varepsilon \left(t\right)=X_t-x_t\left(\vartheta _0\right)$ we have the relation
\begin{align*}
r_\varepsilon \left(t\right)=\int_{0}^{t}\left[ar_\varepsilon
  \left(s\right)+br_\varepsilon \left(\vartheta _0s\right)\right]{\rm d}s+\varepsilon  W_t .
\end{align*}
Therefore 
\begin{align*}
\Ex_{\vartheta _0}\left|r_\varepsilon \left(t\right)\right|^p&\leq
CT^{p-1}\left|a\right|^p\int_{0}^{t}\Ex_{\vartheta _0}\left|r_\varepsilon
\left(s\right)\right|^p{\rm d}s+C\left|b\right|^p\vartheta _0^{-1}\int_{0}^{\vartheta _0t}\Ex_{\vartheta
  _0}\left|r_\varepsilon \left(s\right)\right|^p{\rm d}s\\
&\qquad \qquad +C\varepsilon ^pT^{\frac{p}{2}}\\
&\leq \left[CT^{p-1}\left|a\right|^p+C\left|b\right|^p\vartheta _0^{-1}\right]\int_{0}^{t}\Ex_{\vartheta _0}\left|r_\varepsilon
\left(s\right)\right|^p{\rm d}s+C\varepsilon ^pT^{\frac{p}{2}}
\end{align*}
and
\begin{align*}
R_\varepsilon \left(t\right)\leq \left[CT^{p-1}\left|a\right|^p+C\left|b\right|^p\vartheta _0^{-1}\right]\int_{0}^{t}\Ex_{\vartheta _0}\left|R_\varepsilon
\left(s\right)\right|^p{\rm d}s+C\varepsilon ^pT^{\frac{p}{2}}.
\end{align*}
Now by Gronwall-Bellman lemma we obtain the estimate
\begin{align}
\label{A5}
\Ex_{\vartheta
  _0}\left|X_t-x_t\left(\vartheta _0\right)\right|^p\leq C\,\varepsilon ^p
\end{align}
with the corresponding constant $C>0$.

\section*{Acknowledgments} This research was supported by RSF project no
20-61-47043.

\end{document}